\documentclass[10pt]{article}

\usepackage{amsmath}
\usepackage{amsfonts}
\usepackage{amssymb}
\usepackage{framed}
\usepackage{mathtools}
\usepackage{enumerate}
\usepackage{mathrsfs}
\usepackage[hidelinks]{hyperref}
\usepackage{enumitem}
\usepackage{mathrsfs}
\usepackage{scrextend}
\usepackage{amsthm}
\usepackage{bbm}
\usepackage{thmtools}
\usepackage{thm-restate}
\usepackage{mathabx}
\usepackage[margin=1.5in]{geometry}
\usepackage{rotating}

\usepackage{xcolor}
\hypersetup{
    colorlinks,
    linkcolor={red!50!black},
    citecolor={blue!50!black},
    urlcolor={blue!80!black}
}

\usepackage{tikz-cd}

\theoremstyle{plain}
\newtheorem{thm}{Theorem}[section]
\newtheorem{lem}[thm]{Lemma}
\newtheorem{prop}[thm]{Proposition}
\newtheorem{cor}[thm]{Corollary}

\theoremstyle{definition}
\newtheorem{defn}[thm]{Definition}
\newtheorem{situ}[thm]{Situation}

\theoremstyle{remark}
\newtheorem{rem}[thm]{Remark}

\newtheorem*{notn}{Notation}

\tikzset{
  symbol/.style={
    draw=none,
    every to/.append style={
      edge node={node [sloped, allow upside down, auto=false]{$#1$}}}
  }
}

\newcommand\restr[2]{{
	\left.\kern-\nulldelimiterspace
	#1
	\vphantom{\big|}
	\right|_{#2}
	}}
\newcommand{\an}[1]{#1^{\textrm{an}}}

\newcommand{\ep}{\varepsilon}
\newcommand{\Aut}{\textrm{Aut}}

\newcommand{\ch}[1]{\widecheck{{#1}}}

\newcommand{\NL}{\operatorname{NL}}
\newcommand{\codim}{\operatorname{codim}}

\newcommand{\GL}{\operatorname{GL}}

\newcommand{\SL}{\operatorname{SL}}
\newcommand{\SO}{\operatorname{SO}}
\newcommand{\Sp}{\operatorname{Sp}}

\usepackage{amsmath,calligra,mathrsfs}

\usepackage[cal=boondoxo]{mathalfa}

\DeclareMathOperator{\sheafhom}{\mathcal{H \kern -1pt o \kern -2pt m}}
\DeclareMathOperator{\sheafend}{\mathcal{E \kern -1pt n \kern -2pt d}}
\DeclareMathOperator{\sheafaut}{\mathcal{A \kern -1pt u \kern -2pt t}}

\title{Degrees of Hodge Loci}
\author{David Urbanik}

\begin{document}

\maketitle

\begin{abstract}
We prove asymptotic estimates for the growth in the degree of the Hodge locus in terms of arithmetic properties of the integral vectors that define it. Our methods are general and apply to most variations of Hodge structures for which the Hodge locus is dense. As applications we give asymptotic formulas controlling the degrees of Noether-Lefschetz loci associated to smooth projective hypersurfaces in $\mathbb{P}^3$, and the degrees of subvarieties of the Torelli locus parameterizing Jacobians split up to isogeny.
\end{abstract}

\tableofcontents

\subsection{Preamble}
\label{preamble}

We adopt the convention that all algebraic varieties and schemes are defined over $\mathbb{C}$ unless otherwise stated. The typical exception will be algebraic groups, which will almost always be $\mathbb{Q}$-algebraic. Moreover all Mumford-Tate groups (defined below), will always be understood to be \emph{special} Mumford-Tate groups. For two real sequences $a_{n}$ and $b_{n}$, we write $a_{n} \sim b_{n}$ if $a_{n}/b_{n} \xrightarrow{n \to \infty} 1$ and $a_{n} \asymp b_{n}$ if $a_{n}/b_{n}$ is bounded away from both $0$ and $\pm \infty$ as $n \to \infty$. Finally we write $a_{n} \lesssim b_{n}$ if for any $\ep > 0$ there exists $N$ such that $b_{n} \geq (1 - \ep) a_{n}$ for all $n \geq N$.

\section{Introduction}
\label{intro}

Given a smooth projective family of varieties $f : X \to S$, the cohomology groups $\mathbb{V}_{s} := H^{i}(X_{s}, \mathbb{Z})/\textrm{tor.}$ in degree $i$ carry a natural Hodge structure. Via the Hodge conjecture, rational vectors inside $\mathbb{V}_{s,\mathbb{Q}}$ and its tensor powers are predicted to characterize, up to cohomological equivalence, the algebraic cycles associated to $X_{s}$ and its self-products. A question of much current and classical interest is to understand for which $s$ the fibre $X_{s}$ carries ``more than the expected number'' of algebraic cycles, i.e., more than at a very general point $s \in S(\mathbb{C})$. 

The loci in $S$ where the fibres $X_{s}$ acquire such cycles are examples of Hodge loci (conjecturally, all Hodge loci arise from families of algebraic cycles). Recent work \cite{BKU} \cite{arXiv:2303.16179} has provided a conjectural framework --- and in many cases largely settled --- questions regarding their existence and density. More precisely, following the work in \cite{BKU}, one can give precise conjectures detailing when exactly such loci should exist, what sorts of vectors and tensors should define them, and whether they should be analytically or Zariski dense in $S$. Moreover the work of \cite{arXiv:2211.10592} and \cite{arXiv:2303.16179} allows one to verify such existence and density properties in most cases that arise in practice.

\vspace{0.5em}

However even knowing that such loci are dense in $S$ is in some sense only the beginning of understanding how many such loci there are. A more refined problem is to quantify how such loci grow as one increases the number or type of allowable vectors or tensors which define them. A natural way of counting these loci, dating back to at least \cite{CDK}, is to use a polarizing form $Q : \mathbb{V} \otimes \mathbb{V} \to \mathbb{Z}$, where $\mathbb{V} = R^{i} f_{*} \mathbb{Z} / \textrm{tor}.$ is the natural local system interpolating the fibres $\mathbb{V}_{s}$. Then one can consider the reduced Hodge locus $Z_{n} \subset S$ defined by primitive integral vectors $v$ with $Q(v,v) = n$ for some integer $n$. The locus $Z_{n}$ is algebraic, so after choosing a quasi-projective embedding of $S$ it makes sense to ask for its degree; more generally, it makes sense to ask
\begin{quote}
\textbf{Question:} What are the asymptotics of $\deg Z_{n}$ as $n \to \infty$?
\end{quote}
In this paper we give some general techniques for answering this question.

\subsection{Applications}

To motivate the more technical results that follow, we give some concrete applications of our methods.

\subsubsection{Noether-Lefschetz Loci}
\label{NLexsec}

Let $U_{d} \subset \mathbb{P}^{N}$ be the parameter space of smooth degree $d$ hypersurfaces in $\mathbb{P}^3$ for some integer $d \geq 5$, where $U_{d}$ is the complement in $\mathbb{P}^{N}$ of the discriminant locus $\Delta$. Write $f : X \to S := U_{d}$ for the universal family. Using the natural map $i : X \to \mathbb{P}^{3}$ we obtain a family $\mathcal{L} := i^{*} \mathcal{O}(1)$ of line bundles in the fibres of $f$. We define the primitive subsystem
\[ \mathbb{V} := \ker \left[ R^{2} f_{*} \mathbb{Z} \xrightarrow{(-) \cup c^{1}(\mathcal{L})}  R^{4} f_{*} \mathbb{Z} \right] \]
where $c^{1}(\mathcal{L})$ is the global section of $R^{2} f_{*} \mathbb{Z}$ coming from the relative Chern class of $\mathcal{L}$. The local system $\mathbb{V}$ is torsion-free, and cup product induces a non-degenerate symmetric pairing $Q : \mathbb{V} \otimes \mathbb{V} \to \mathbb{Z}$ which makes, as a consequence of Poincar\'e duality, each fibre $(\mathbb{V}_{s}, Q_{s})$ into a unimodular lattice. We write $(r,s)$ for the signature of this lattice, which is related to the Hodge numbers $(h^{2,0}, h^{1,1}, h^{0,2})$ of $\mathbb{V}$ by $r = 2 h^{2,0} = 2 h^{0,2}$ and $s = h^{1,1}$.

The Noether-Lefschetz locus of $\mathbb{V}$ is constructed as follows. For each point $s \in S(\mathbb{C})$ and each integral vector $v \in \mathbb{V}_{s}$, there is a locus $Z(v) \subset S$ consisting of all $s' \in S(\mathbb{C})$ where some flat translate of $v$ in $\mathbb{V}_{s'}$ is Hodge. By \cite{CDK} the locus $Z(v)$ is algebraic, and we will choose to equip it with its underlying reduced structure. Moreover \cite{CDK} also shows that for $\ell \in \mathbb{Z}$ the locus $\bigcup_{Q(v,v) = \ell} Z(v)$ is algebraic.

\begin{defn}
We say a vector $v$ in a $\mathbb{Z}$-lattice $L$ is primitive if $v$ has non-zero image in $L_{\mathbb{F}_{p}}$ for each $p$.
\end{defn}

Working with the union over all $v$ such that $Q(v,v) = \ell$ results in overcounting: if $\lambda \in \mathbb{Z}$ then we have $Z(v) = Z(\lambda v)$. To remedy this, we define a modified scalar-invariant quantity better suited for counting Hodge loci which agrees with $Q(v,v)$ whenever $v$ is primitive. 

\begin{defn}
\label{Qbardef}
For $v$ a non-zero vector in some fibre of $\mathbb{V}_{\mathbb{Q}}$, we define $a(v) \in \mathbb{Q}^{\times}$ to be the unique positive rational scalar for which $a(v) v$ is a primitive integral vector in the corresponding fibre of $\mathbb{V}$. We moreover define
\[ \overline{Q} : \mathbb{V}_{\mathbb{Q}} \setminus \{ 0 \} \to \mathbb{Z}, \hspace{2em} \overline{Q}(v) = a(v)^2 Q(v,v) = Q(a(v) v, a(v) v) \]
which is a map of flat bundles over $S$. We write $u(v)$ and $\nu(v)$ for the unique integers such that $\overline{Q}(v) = u(v) \nu(v)^2$ with $\nu(v)$ positive and $u(v)$ square-free. 
\end{defn}

\noindent It is clear that $\overline{Q}(\lambda v) = \overline{Q}(v)$ for each $\lambda \in \mathbb{Q}^{\times}$, and that $u(v) = \overline{Q}(v) = Q(v,v)$ modulo scaling by $(\mathbb{Q}^{\times})^2$. We may then define
\begin{align*}
\NL_{\nu,u} &:= \left\{ s \in S(\mathbb{C}) : \begin{array}{c}\textrm{ there exists a Hodge vector }v \in \mathbb{V}_{s} \\ \textrm{ with }u(v) = u \textrm{ mod }(\mathbb{Q}^{\times})^2, \, \, \nu(v) | \nu \end{array} \right\} \\
&= \left\{ s \in S(\mathbb{C}) : \begin{array}{c}\textrm{ there exists a Hodge vector }v \in \mathbb{V}_{s} \\ \textrm{ with }Q(v,v) = u \textrm{ mod }(\mathbb{Q}^{\times})^2, \, \, \nu(v) | \nu \end{array} \right\} .
\end{align*}
For a closed subvariety of $U_{d}$ we define its degree to be the degree of its closure in $\mathbb{P}^{N}$. We set $S = U_{d}$ and $m = r+s$.

\begin{thm}[\textbf{Upper bound}]
\label{NLlocidegthmupper}
For each fixed square-free $u \in \mathbb{Z}_{> 0}$ we have
\[ \deg(\NL_{\nu,u}) \lesssim c_{u} \nu^{m+2h^{0,2}-2} \]
as $\nu \to \infty$, where $c_{u} > 0$ is a real constant depending on $u$.
\end{thm}

\vspace{0.5em}

\begin{rem}
The quantity $m$ is $b_{2} - 1$, where $b_{2}$ is the second Betti number of a fibre of $f$. It is known that $b_{2} = d^3 - 4d^2 + 6d - 2$, and by \cite[Ch 17.3]{zbMATH05982705} one has that $h^{0,2} = (d-1)(d-2)(d-3)/6$. Combining these estimates one easily expresses $m + 2h^{0,2} - 2$ as a polynomial in $d$.
\end{rem}

\begin{defn}
The tensorial Hodge locus of $\mathbb{V}$ is the collection of all Hodge loci (i.e., Noether-Lefschetz loci) associated to the variations $\mathbb{V}^{\otimes k}$ for all $k \geq 0$.
\end{defn}

\vspace{0.5em}

We say a subvariety of $S$ is Hodge-generic if it is not contained in the tensorial Hodge locus. In what follows the term ``sufficiently general hyperplane section'' is understood as in \autoref{suffgenhypplanedef}. For the use of the term ``period dimension'' in the following statement, see \autoref{perdimdef}.

\begin{thm}[\textbf{Lower Bound}] Fix the data of:
\label{NLlocidegthmlower}
\begin{itemize}
\item[-] a locally closed irreducible Hodge-generic subvariety $T \subset S = U_{d}$, with period dimension $\geq h^{2,0} + 1$;
\item[-] a relatively compact open neighbourhood $B \subset T$; and
\item[-] a sufficiently general hyperplane section $L \cap B$ of $B$ with $\codim_{S} L = \dim T - h^{2,0}$.
\end{itemize}
Then for each fixed square-free $u \in \mathbb{Z}_{> 0}$, and for $\nu$ with $\nu \in \mathbb{Z}_{> 0} \setminus 2\mathbb{Z}$, one has 
\begin{align*}
c_{u} \nu^{m-2} &\lesssim | \NL_{\nu,u} \cap L \cap B | 
\end{align*}
as $\nu \to \infty$, where $c_{u} \in \mathbb{R}_{> 0}$ is a constant depending only on $u, B$ and $L$. 
\end{thm}

In the special case where $S = T$, one has

\begin{thm}[\textbf{Lower Bound}]
\label{pureNLlowerbound}
For each fixed square-free $u \in \mathbb{Z}_{> 0}$, and for $\nu$ with $\nu \in \mathbb{Z}_{> 0} \setminus 2\mathbb{Z}$, one has 
\begin{align*}
c_{u} \nu^{m-2} &\lesssim \deg \NL_{\nu,u}
\end{align*}
as $\nu \to \infty$, where $c_{u} \in \mathbb{R}_{> 0}$ is a real constant.
\end{thm}

\begin{rem}
The lower bound given in \autoref{NLlocidegthmlower} can be made optimal at the cost of a more complicated algebraic expression in the prime factorization of $\nu$. We give the precise calculation in \S\ref{volcompsubsubsec}. With the optimal formula appearing there, one can replace the symbol $\lesssim$ with $\asymp$ as long as $B$ does not intersect a certain finite union of tensorial Hodge loci of $T$. 
\end{rem}

\vspace{0.5em}

It is a consequence of the theory developed in \cite{BKU} (c.f. \cite[Thm. 6]{nondensityexp}) that each component of $\NL_{\nu,u} \cap T$ has codimension exactly $h^{2,0}$ in $T$ outside of a closed algebraic locus in $T$, so away from this locus the ``sufficiently general'' hypothesis in \autoref{NLlocidegthmlower} ensures that the cardinality $| \NL_{\nu,u} \cap L \cap B |$ is finite. The left-hand side of the inequality in \autoref{NLlocidegthmlower} diverges as a function of $\nu$, so one obtains infinitely many Hodge loci defined by vectors $v$ with $u(v) = u$ in any neighbourhood of $T$. In fact, away from the atypical Hodge locus, one can obtain exact asymptotics in terms of $\nu$ for the size of the Hodge locus in the prescribed region.

The period dimension hypothesis of $\geq h^{2,0} + 1$ on the other hand implies that each component of $\NL_{\nu, u} \cap T$ has positive period dimension; this will be used to verify that a general component of $\NL_{\nu, u} \cap T$ carries a unique global $1$-dimensional subspace of Hodge vectors. For the importance of carrying out such a verification see our discussion in \S\ref{redtoptcountsec} and \S\ref{compwithotherworksec}.

\subsubsection{Split Jacobians}
\label{splitjacexsec}

In this case we start with the universal family $f : X \to S := \mathcal{M}_{g}$ of genus $g$ curves (strictly speaking, we make sense of this using the language of Deligne-Mumford stacks, but the results would be no different were we to rigidify the moduli problem by adding marked points). We fix a projective compactification $\overline{\mathcal{M}_{g}}$ of (the coarse space of) $\mathcal{M}_{g}$ and define the degree of an algebraic locus in $\mathcal{M}_{g}$ to be the degree of its closure in $\overline{\mathcal{M}_{g}}$. 

\begin{defn}
For a complex abelian variety $A$, we say that $A$ has a $k$-factor of degree $\nu$ if there is an isogeny $\varphi : A \to A' \times A''$ such that $\deg \varphi = \nu$ and $\dim A' = k$.
\end{defn}
For a curve $C$, we write $J(C)$ for the associated Jaocbian variety. As an analogue of the Noether-Lefschetz loci considered above, we consider
\begin{align*}
\textrm{SP}_{k,\nu} := \{ s \in S(\mathbb{C}) : J(X_{s}) \textrm{ has an }k\textrm{-factor of degree dividing }\nu \} . 
\end{align*}

\begin{thm}[\textbf{Upper Bound}]
\label{abidemthmupper}
For each integer $k$ such that $1 \leq k \leq g/2$ we have
\[ \deg(\textrm{SP}_{k,\nu}) \lesssim c_{k} \nu^{8 k(g-k)} \]
as $\nu \to \infty$, where $c_{k} > 0$ is a positive real constant.
\end{thm}

We define the atypical Hodge locus and the notion of a ``sufficiently general'' hyperplane as in the previous section except with $\mathbb{V} = R^{1} f_{*} \mathbb{Z}$. We consider the pair $(\Sp_{g}, \mathbb{H}_{g})$ consisting of the symplectic group of a $2g$-dimensional integral symplectic space $(V,Q)$, and the Siegel upper half space, respectively. This pair can be viewed as the ``generic Hodge datum'' of the variation of Hodge structure $\mathbb{V}$ (we review this notion below). The complex manifold $\mathbb{H}_{g}$ has natural closed submanifolds $\mathbb{H}_{k} \times \mathbb{H}_{g-k}$ whose images in $\mathcal{A}_{g}(\mathbb{C}) = \Sp_{g}(\mathbb{Z}) \backslash \mathbb{H}_{g}$ consist exactly of those closed algebraic loci parameterizing non-simple principally polarized abelian varieties. Each such submanifold corresponds to a decomposition $V_{\mathbb{Q}} \simeq V_{1} \oplus V_{2}$, compatible with the symplectic form, and the associated symplectic idempotents $e_{1}$ and $e_{2}$ belong to only finitely many $\Sp_{g}(\mathbb{R})$ orbits as the submanifolds $\mathbb{H}_{k} \times \mathbb{H}_{g-k} \hookrightarrow \mathbb{H}_{g}$ vary. We let $O_{k}$ be the $\Sp_{g}$-orbit of such an idempotent $e_{k}$ whose image has dimension $2k$. We view $O_{k}$ as a $\mathbb{Q}$-algebraic subvariety of the affine space associated to the free $\mathbb{Z}$-module $\textrm{End}(V)$. 

\begin{thm}[\textbf{Lower Bound}]
\label{abidemthm}
Let $1 \leq k \leq g-1$ and $g \geq 2$. Fix the data of: 
\begin{itemize}
\item[-] a locally closed irreducible Hodge generic subvariety $T \subset S$ of dimension at least $k(g-k) + 1$;
\item[-] a relatively compact open neighbourhood $B \subset T$; and
\item[-] a sufficiently general hyperplane section $L \cap B$ of $B$ with $\codim_{S} L = \dim T - k(g-k)$.
\end{itemize}
Then with $O = O_{k}$ one has 
\begin{equation}
\label{ablowboundeq}
 c_{k} \prod_{p} \mu_{p,O}\left( O(\mathbb{Q}_{p}) \cap \frac{1}{\nu} \textrm{End}(V)_{\mathbb{Z}_{p}} \right) \lesssim |\textrm{SP}_{k,\nu^{2g}} \cap L \cap B|
\end{equation}
as $\nu \to \infty$, where
\begin{itemize}
\item[-] $c_{k}$ is a positive real constant depending only on $k, B$ and $L$; and
\item[-] $\mu_{p,O}$ is the unique $\Sp_{g}(\mathbb{Q}_{p})$-invariant measure on $O(\mathbb{Q}_{p})$ compatible with the Tamagawa measures on $\Sp_{g}$ and the stabilizer in $\Sp_{g}$ of $e_{k}$.
\end{itemize}
\end{thm}

\begin{rem}
In the proof of \autoref{abidemthm}, it will be more natural to work with a certain locus $\textrm{SP}'_{k,\nu}$ defined Hodge-theoretically. We have $\textrm{SP}'_{k,\nu} \subset \textrm{SP}_{k,\nu^{2g}}$. Moreover if $B$ does not intersect a certain finite union of tensorial Hodge loci of $T$, then the statement of \autoref{abidemthm} holds with $\textrm{SP}_{k,\nu^{2g}}$ replaced by $\textrm{SP}'_{k,\nu}$ and the symbol $\lesssim$ replaced by $\asymp$. 
\end{rem}

\begin{rem}
One can show that the left-hand side of (\ref{ablowboundeq}) diverges as $\nu \to \infty$. 
\end{rem}

Let us take some time to unpack the statement. The group $\Sp_{g}$ is semisimple, simply-connected, and unimodular. By \cite[pg.26]{zbMATH03192559} one also knows that the locally compact topological group $\Sp_{g}(\mathbb{A})$ is unimodular, where we write $\mathbb{A}$ for the adele ring. The semisimplicity and simply-connected properties imply that $\Sp_{g}(\mathbb{Q})$ is a lattice in $\mathbf{H}_{S}(\mathbb{A})$ \cite[Thm. 1]{zbMATH03194969} \cite{zbMATH03194968}, and this lets us put a unique Haar measure $\mu$, the Tamagawa measure, on $\Sp_{g}(\mathbb{A})$ such that $\Sp_{g}(\mathbb{A})/\Sp_{g}(\mathbb{Q})$ has unit volume. The measure $\mu$ splits as a product $\mu_{\infty} \times \prod_{p} \mu_{p}$ of measures defined at each place of $\mathbb{Q}$.

The stabilizer in $\Sp_{g}$ of a Hodge-theoretic idempotent $e_{k}$ is abstractly isomorphic to a product $H = \textrm{Sp}_{k} \times \textrm{Sp}_{g-k}$ of symplectic groups, and the same reasoning applies, mutatis mutandis, to $H$. Moreover by \cite[A.1.2, A.5.1]{zbMATH05902967} the space $(\Sp_{g}/H)(\mathbb{A})$ consists of finitely many $\Sp_{g}(\mathbb{A})$ orbits; let us fix such an orbit $(\Sp_{g}/H)(\mathbb{A})^{\circ}$. Using \cite[Thm. 5.2.1, 5.2.2]{zbMATH06960207} there is a unique invariant measure on this quotient determined by the measures $\mu$ on $\Sp_{g}(\mathbb{A})$ and $H(\mathbb{A})$ which again splits as a product of measures associated to each place of $\mathbb{Q}$. Using the natural identification $\Sp_{g}/H \simeq O = O_{k}$ this induces in a natural way a measure on a component $O(\mathbb{A})^{\circ} \subset O(\mathbb{A})$ and on $O(\mathbb{Q}_{p})$ for each $p$ (note here that $O(\mathbb{Q}_{p})$ consists of a single $\Sp_{g}(\mathbb{Q}_{p})$ orbit by \cite[A.5.3]{zbMATH05902967}), and it is these measures we use in the statement of \autoref{abidemthm}.

\subsubsection{Comparing The Examples}

To understand the relationship between our examples, let us explain how one would actually go about computing the asymptotic expression in \autoref{abidemthm}. The Haar measure on the group $\Sp_{g}(\mathbb{A})$ can also be described using a $\Sp_{g}$-invariant $\mathbb{Q}$-algebraic form $\omega$ on $\Sp_{g}$ of top degree, in such a way so that the measures at each place of $\mathbb{Q}$ are obtained by integrating against $\omega$. Such a form is called a \emph{gauge} form. It turns out that the measures associated to $H$ and the quotient $\Sp_{g}/H$ are also associated to gauge forms, and so the quantities appearing in the product in \autoref{abidemthm} can be computed assuming one is able to compute the relevant volume integrals over the $p$-adic manifolds $O(\mathbb{Q}_{p})$.

In the Noether-Lefschetz locus case, essentially the same reasoning applies. In that case, the relevant group is $\SO(r,s)$ instead of $\Sp_{g}$. In particular, up to the scaling action of $\mathbb{Q}^{\times}$, the orbits of the group $\SO(r,s)(\mathbb{R})$ on fibres of $\mathbb{V}$ are actually classified by the invariant $u$, with two rational vectors $v, v'$ lying in the same $\mathbb{Q}^{\times} \SO(r,s)(\mathbb{R})$-orbit if and only if $u(v) = u(v')$. (This follows from the fact that $u(v) = Q(v,v) = \overline{Q}(v)$ modulo $(\mathbb{Q}^{\times})^2$, and then by using that $u(v)$ is the unique integral square-free representative of the coset $u(v) (\mathbb{Q}^{\times})^2$.) The product lower bound in \autoref{NLlocidegthmlower} is in fact derived from an analogous product of adelic volumes coming from a natural adelic measure on a corresponding orbit of $\SO(r,s)$. The difference being that, in this case, we have taken the time to actually compute the relevant $p$-adic volume integrals, and the factor at $p$ gives a contribution lower bounded by $p^{(m-2) e_{p}(\nu)}$, where $e_{p}(\nu)$ is the exponent of $p$ in the prime factorization of $\nu$.

We expect similar computations can be done to give an explicit form to the left-hand side of \autoref{abidemthm}, but this task is beyond the scope of this paper.

\subsection{The General Case}

\subsubsection{Variational Setup}

To facilitate the above applications, we now describe some general theorems applying to an arbitrary integral variation of Hodge structure $\mathbb{V}$ with polarization $Q$ on a smooth quasi-projective complex variety $S$. This will require some additional setup. We will simplify the situation by studying the asymptotics only of Hodge loci associated to vectors of $\mathbb{V}$; the case of loci associated to general Hodge tensors reduces to this by replacing $\mathbb{V}$ with some tensor construction. We may also assume the weight of $\mathbb{V}$ is zero, Tate-twisting if necessary. In particular we have a Hodge filtration $F^{\bullet}$ on $\mathcal{H} := \mathbb{V} \otimes_{\mathbb{Z}} \mathcal{O}_{\an{S}}$ and the space of Hodge vectors above $s$ is the subspace $F^{0}_{s} \cap \mathbb{V}_{\mathbb{Q}, s}$. Given a vector $v$ in some fibre of $\mathbb{V}_{\mathbb{Q}}$, we write $Z(v)$ for the corresponding (reduced) Hodge locus. We define $\overline{Q}$, $u$ and $\nu$ as in \autoref{Qbardef}. 

For each point $s \in S(\mathbb{C})$ we write $h_{s}$ for the Hodge structure on $\mathbb{V}_{s}$, which we view as a map $h_{s} : \mathbb{U} \to \Aut(\mathbb{V}_{s}, Q_{s})_{\mathbb{R}}$. Here $\mathbb{U}$ is the ``circle group'', a real-algebraic subgroup of the Deligne torus $\mathbb{S} = \textrm{Res}_{\mathbb{C}/\mathbb{R}} \mathbb{G}_{m,\mathbb{C}}$ constructed as in \cite[I.A]{GGK}. The Mumford-Tate group $\mathbf{G}_{s}$ at $s$ is the $\mathbb{Q}$-Zariski closure of $h_{s}(\mathbb{U})$, and the Mumford-Tate domain $D_{s}$ at $s$ is the $\mathbf{G}_{s}(\mathbb{R})$-conjugacy class of the morphism $h_{s} : \mathbb{U} \to \mathbf{G}_{s,\mathbb{R}}$. We call the pair $(\mathbf{G}_{s}, D_{s})$ the Hodge datum at $s$. The generic Hodge datum $(\mathbf{G}_{S}, D_{S})$ is the is the abstract Hodge datum which is isomorphic to $(\mathbf{G}_{s}, D_{s})$ at a point $s \in S(\mathbb{C})$ outside the Hodge locus of all the tensor powers $\mathbb{V}^{\otimes k}$ for $k \geq 1$. The flat structure of the local system locally identifies the realizations $(\mathbf{G}_{S,s}, D_{S,s})$ of the generic Hodge datum and allows one to view the pair $(\mathbf{G}_{S}, D_{S})$ as an abstract object with a realization in each fibre of $\mathbb{V}$. 

\begin{defn}
For a local system $\mathbb{L}$ on a smooth complex algebraic variety $S$, the algebraic monodromy group $\mathbf{H}_{\mathbb{L},S,s}$ is the identity component of the Zariski closure of the map
\[ \pi_{1}(S, s) \to \GL(\mathbb{L}_{s})  \]
where $s \in S(\mathbb{C})$ is a point. Usually either some base point $s$, the local system $\mathbb{L}$, or both, are understood, in which case we write $\mathbf{H}_{\mathbb{L},S}, \mathbf{H}_{S,s}$ or $\mathbf{H}_{S}$, as appropriate. 
\end{defn}

Write $\mathbf{H}_{S}$ for the algebraic monodromy group of $\mathbb{V}$. By \cite[\S5]{Andre1992}, $\mathbf{H}_{S}$ is a normal subgroup of the derived subgroup $\mathbf{D}_{S}$ of $\mathbf{G}_{S}$. We note that, by \cite[Lem. 2]{Andre1992}, the group $\mathbf{G}_{S}$ is always at least reductive, which implies that $\mathbf{H}_{S}$ is always semisimple. 

\subsubsection{Period Maps}

An important part of formulating a generalization of the examples we have discussed is deciding when a dense collection of Hodge loci should exist in $S$; this is not automatic, since for variations of higher level (in the sense of \cite[\S4.6]{BKU}) one expects by \cite[Conj. 3.5]{BKU} that the Hodge locus should lie in a strict algebraic subvariety of $S$. Recently a general criterion for proving the existence of such loci was given first in \cite{arXiv:2211.10592} in an abstract setting, and sharpened in the Hodge-theoretic setting in \cite{arXiv:2303.16179}. One of our main technical results will in some sense be a further strengthening of \cite{arXiv:2303.16179}, so we start by adopting some of the setup. Because the statement of our main results can be deduced after replacing $S$ with a finite cover $g : S' \to S$ and $\mathbb{V}$ with $\mathbb{V}' = g^{*} \mathbb{V}$, we will freely make such changes.

\paragraph{Associated Period Maps:} We may replace $S$ with a finite \'etale covering such that $\mathbb{V}$ induces a period map $\varphi : S \to \Gamma \backslash D_{S}$, where $\Gamma \subset \mathbf{G}_{S}(\mathbb{Q})$ is a fixed neat arithmetic lattice containing the image of monodromy. Moreover by applying \cite[III.A.1]{GGK} we may assume that we have a factorization $\varphi = \varphi_{1} \times \cdots \times \varphi_{n}$, such that
\begin{itemize}
\item[(i)] the factorization of $\varphi$ is induced by a factorization $\Gamma \backslash D_{S} = \Gamma_{1} \backslash D_{1} \times \cdots \times \Gamma_{n} \backslash D_{n}$ with $\Gamma = \Gamma_{1} \times \cdots \times \Gamma_{n}$; and
\item[(ii)] the factorizations in (i) are induced by the almost-direct product factorization $\mathbf{D}_{S} = \mathbf{H}_{1} \cdots \mathbf{H}_{n}$, where we note that the center of $\mathbf{G}_{S}$ acts trivially on $D_{S} = D_{1} \times \cdots \times D_{n}$.
\end{itemize}
Given a subset $I \subset \{ 1, \hdots, n \}$, we write $p_{I} : \Gamma \backslash D_{S} \to \Gamma_{I} \backslash D_{I}$ where $\Gamma_{I} = \prod_{i \in I} \Gamma_{i}$ and $D_{I} = \prod_{i \in I} D_{i}$. The map $\varphi_{I}$ is defined similarly. 

\begin{defn}
\label{perdimdef}
For a subvariety $T \subset S$, its period dimension is the dimension of $\varphi(T)$. 
\end{defn}

\paragraph{Likeliness:} Given the generic Hodge datum $(\mathbf{G}_{S}, D_{S})$, a Hodge sub-datum is a pair $(\mathbf{M}, D_{M})$ consisting of a subgroup $\mathbf{M} \subset \mathbf{G}_{S}$ which is the Mumford-Tate group of some Hodge structure $h \in D_{S}$ and the associated orbit $D_{M} = \mathbf{M}(\mathbb{R}) \cdot h$ under the conjugation action. The Hodge datum of a point $h \in D_{S}$ is the pair $(\mathbf{M}_{h}, \mathbf{M}_{h}(\mathbb{R}) \cdot h)$, where $\mathbf{M}_{h}$ is the Mumford-Tate group of $h$.

\begin{defn}
We say the Hodge sub-datum $(\mathbf{M}, D_{M}) \subset (\mathbf{G}_{S}, D_{S})$ is $\mathbb{V}$-likely if for every non-empty set $I \subset \{ 1, \hdots, n \}$ the inequality $\dim \pi_{I}(D_{M}) + \dim \varphi_{I}(\an{S}) - \dim D_{I} \geq 0$ holds, and we say it is \emph{strongly} $\mathbb{V}$-likely if these inequalities are strict.
\end{defn}

\paragraph{Subdata of Definition:} A component $Z \subset S$ of the (tensorial) Hodge locus of $S$ can be defined by many different Hodge subdata. This means the following:

\begin{defn}
\label{Hodgedatadef}
We say a Hodge subdatum $(\mathbf{M}, D_{M}) \subset (\mathbf{G}_{S}, D_{S})$ \emph{defines} a closed algebraic locus $Z \subset S$ if $Z$ is an irreducible component of $\varphi^{-1}(\pi(D_{M}))$.
\end{defn}

\noindent Note that $\mathbf{M}$ need not be the Mumford-Tate group of $Z$ in (\ref{Hodgedatadef}). A Hodge locus $Z$ is always defined by a Hodge subdatum $(\mathbf{G}_{Z}, D_{Z}) \subset (\mathbf{G}_{S}, D_{S})$ which arises from the pair $(\mathbf{G}_{s}, D_{s})$ at a very general point $s \in Z(\mathbb{C})$, but it is often useful to view $Z$ as defined by a potentially larger subdatum.

\subsubsection{Equivalence Classes of Hodge Subdata}

The group $\mathbf{G}_{S}(\mathbb{R})$ acts on Hodge subdata $(\mathbf{M}, D_{M}) \subset (\mathbf{G}_{S}, D_{S})$ by conjugation, with $(\mathbf{M}, D_{M})$ being equivalent to $(\mathbf{M}', D_{M'})$ if there is $g \in \mathbf{G}_{S}(\mathbb{R})$ such that $\mathbf{M}' = g \mathbf{M} g^{-1}$ and $D_{M'} = g D_{M}$. We say Hodge data that are equivalent in this way are conjugation-equivalent. There are finitely many conjugation equivalence classes of Hodge subdata, a fact which is explained in \cite[\S4.4]{arXiv:2303.16179}.

\begin{defn}
\label{cosetequiv}
We say that two Hodge subdata $(\mathbf{M}, D_{M}), (\mathbf{M}', D_{M'}) \subset (\mathbf{G}_{S}, D_{S})$ are $\mathbb{Q}$-coset equivalent if there exists $g \in \mathbf{G}_{S}(\mathbb{R})$ such that $(g \mathbf{M}' g^{-1}, g D_{M'}) = (\mathbf{M}, D_{M})$, and the coset $g \mathbf{M}'$ defines a rational point in $(\mathbf{G}_{S}/\mathbf{M}')(\mathbb{Q})$. 
\end{defn}

\begin{lem}
$\mathbb{Q}$-coset equivalence is an equivalence relation on Hodge subdata.
\end{lem}

\begin{proof}
We observe that the notion of being $\mathbb{Q}$-coset equivalent is symmetric: if $g \mathbf{M}'$ defines a rational point in $\mathbf{G}_{S}/\mathbf{M}'$, then $g \mathbf{M}' = \mathbf{M} g$ is a $\mathbb{Q}$-algebraic subvariety of $\mathbf{G}_{S}$. Taking inverses, this implies that $g^{-1} \mathbf{M}$ is a $\mathbb{Q}$-algebraic subvariety of $\mathbf{G}_{S}$, and hence corresponds to a point of $(\mathbf{G}_{S}/\mathbf{M})(\mathbb{Q})$. 

For transitivity, we consider three tuples $(\mathbf{M}, D_{M}), (\mathbf{M}', D_{M'})$ and $(\mathbf{M}'', D_{M''})$ related by $\mathbb{Q}$-cosets $g \mathbf{M}$ and $g' \mathbf{M}'$. Then for any automorphism $\sigma \in \textrm{Aut}(\mathbb{C}/\mathbb{Q})$ we have
\begin{align*}
(g' g \mathbf{M})^{\sigma} &= g'^{\sigma} g \mathbf{M} \\
&= g'^{\sigma} \mathbf{M}' g \\
&= g' \mathbf{M}' g = g' g \mathbf{M} .
\end{align*}
\end{proof}

\noindent A given real conjugacy class of Hodge subdata, when further partitioned by $\mathbb{Q}$-coset equivalence, can either result in countably-infinitely many subequivalence classes or finitely many. The situation with Noether-Lefschetz loci in \S\ref{NLexsec} is an example of the first case, where the subequivalence classes are indexed by the invariant $u$, and the second case is exhibited by the split Jacobian loci discussed in \S\ref{splitjacexsec}. 

A useful converse result is the following:

\begin{prop}
\label{Qpointgivesdatum}
Let $(\mathbf{M}, D_{M}) \subset (\mathbf{G}_{S}, D_{S})$ be a Hodge subdatum, let $g \in \mathbf{G}_{S}(\mathbb{R})$, and suppose that $g \mathbf{M} \in (\mathbf{G}_{S}/\mathbf{M})(\mathbb{Q})$. Then $(g \mathbf{M} g^{-1}, g D_{M})$ is a Hodge subdatum of $(\mathbf{G}_{S}, D_{S})$. 
\end{prop}

\begin{proof}
See \autoref{Htransconstrlem}. 
\end{proof}

\noindent In many situations in practice, this means that counting $\mathbb{Q}$-coset equivalent Hodge subdata is the same as counting rational points in $\mathbf{G}_{S}/\mathbf{M}$. 

\subsubsection{From Hodge Loci to Point Counting}

We are now ready to describe our main technical results. 

\begin{situ}
\label{mainthmsit}
Let $\mathbb{V}$ be a polarized integral variation of Hodge structure on $S$ satisfying $\mathbf{H}_{S} = \mathbf{D}_{S}$. Suppose that $(\mathbf{M}, D_{M}) \subset (\mathbf{G}_{S}, D_{S})$ is a Hodge subdatum, and that under the natural Hodge representation $\rho : \mathbf{G}_{S} \to \GL(\mathbb{V}_{s})$
\begin{itemize}
\item[-] the group $\mathbf{M}$ is identified with the stabilizer of some $v \in \mathbb{V}_{s}$;
\item[-] one has $(\mathbf{G}_{s}, D_{s}) \subset (\mathbf{M}, D_{M})$, where $(\mathbf{G}_{s}, D_{s})$ is the Hodge datum at $s$.
\end{itemize}
Let $O = \rho(\mathbf{G}_{S}) \cdot v$ be the orbit of $v$ and write $O(\mathbb{R})^{\circ} \subset O(\mathbb{R})$ for the component containing $v$, and set $O(\mathbb{Q})^{\circ} = O(\mathbb{Q}) \cap O(\mathbb{R})^{\circ}$. Fix a projective compactification $S \subset \overline{S}$ and define the degree of $Z \subset S$ to be the degree of its closure.
\end{situ}

Note that the orbit $O$ is preserved by the action of $\Gamma_{S}$, and so has, for each $s \in S(\mathbb{C})$, a well-defined realization $O_{s}$ as a subvariety of the affine space associated to $\mathbb{V}_{s}$. Likewise it makes sense to consider $O(\mathbb{Q}) \cap \frac{1}{\nu} \mathbb{V}$ in each fibre of $\mathbb{V}$: those points $v' \in O(\mathbb{Q})$ with the property that $\nu v'$ is integral. (Here we are using the notation ``$\mathbb{V}$'' in the same spirit as with $\mathbf{G}_{S}$ and $D_{S}$ above: it is an abstract object with a realization above each point $s \in S(\mathbb{C})$.) Given a vector $v' \in O(\mathbb{Q})$ we obtain, using the fixed isomorphism $\mathbf{G}_{S}/\mathbf{M} \simeq O$, a Hodge subdatum $(\mathbf{M}', D_{M'})$ defined as in \autoref{Qpointgivesdatum}. We then set $Z(v') = \varphi^{-1}(\pi(D_{M'}))$ to be the corresponding Hodge locus. 

\begin{thm}[\textbf{Upper Bound}]
\label{mainthmupbound}
Work in the situation of \autoref{mainthmsit}. Then there exists:
\begin{itemize}
\item[-] finitely many Hodge-theoretic Siegel sets $\mathfrak{G}_{1}, \hdots, \mathfrak{G}_{\ell}$ for $\mathbf{G}_{S}$ associated to Hodge structures $h_{1}, \hdots, h_{\ell} \in D_{S}$, such that the Hodge datum of each $h_{i}$ is contained in a Hodge datum $\mathbb{Q}$-coset equivalent to $(\mathbf{M}, D_{M})$; and
\item[-] Hodge vectors $v_{1}, \hdots, v_{\ell} \in O(\mathbb{Q})^{\circ}$, with $v_{i}$ a Hodge vector for $h_{i}$;
\end{itemize}
such that we have
\begin{equation}
\label{degpointcorres}
\deg \left[ \bigcup_{\substack{v' \in \mathcal{O}(\mathbb{Q})^{\circ} \\ \nu(v') | \nu}} Z(v') \right] \lesssim c^{+} \sum_{i=1}^{\ell} \left| \mathfrak{G} \cdot v_{i} \cap \frac{1}{\nu} \mathbb{V} \right| ,
\end{equation}
as $\nu \to \infty$, where $c^{+} > 0$ is a real constant depending on $S, \mathbb{V}$ and $Q$. 
\end{thm}

\begin{notn}
For each fixed non-negative integer $d$ and analytic variety $Z$, we write $Z_{d}$ for the union of $d$-dimensional components of $Z$.
\end{notn}

\begin{thm}[\textbf{Lower Bound}]
\label{mainthmlowbound}
Work in the situation of \autoref{mainthmsit}, assume that $(\mathbf{M}, D_{M})$ is strongly $\mathbb{V}$-likely, and that $\mathbf{M}$ has finite index in its normalizer in $\mathbf{G}_{S}$. Let $B \subset S(\mathbb{C})$ be any open neighbourhood. Then there exists a sufficiently general hyperplane section $L \subset S$ which intersects $B$ and has codimension $\dim S - (\dim D_{S} - \dim D_{M})$, and an open subset $U \subset O(\mathbb{R})^{\circ}$ such that
\begin{equation}
\label{generalowerboundeq}
c^{-} \left| U \cap \frac{1}{\nu} \mathbb{V} \right| \lesssim \sum_{\substack{v' \in \mathcal{O}(\mathbb{Q})^{\circ} \\ \nu(v') | \nu }} |(Z(v') \cap L \cap B)_{0}|
\end{equation}
as $\nu \to \infty$, where $c^{-} > 0$ is a real constant depending on $S, \mathbb{V}$ and $Q$. 
\end{thm}

\vspace{0.5em}

The basic summary of the the two theorems is that Hodge loci defined by vectors $v'$ with $\nu = \nu(v')$ are roughly in correspondence with rational vectors in certain subsets of $O(\mathbb{R})^{\circ}$ with denominators dividing $\nu$. These subsets are at worst as large as a Siegel set orbit, and in the strongly $\mathbb{V}$-likely case contain open neighbourhoods corresponding to any chosen open subset of $S(\mathbb{C})$. The Siegel sets $\mathfrak{G}_{i}$ referred to in \autoref{mainthmupbound} are ``Hodge-theoretic'' in the sense that their associated maximal compact subgroups $K_{i}$ are naturally constructed from the Hodge metric at a point $h_{i} \in D_{S}$; we give a precise description in \S\ref{siegesetsec}. Moreover we act with these sets on Hodge vectors $v_{i}$ for the associated $h_{i}$. 

It turns out that such Siegel-set orbits constructed from Hodge-theoretic data have highly constrained geometry: they are in general non-compact, and have infinite volume for the natural invariant measure on $O(\mathbb{R})^{\circ}$, but nevertheless satisfy the property that, for any positive integer $\nu$, they contain only finitely many rational points with denominators of size at most $\nu$. This means in particular that counting points in such orbits will make sense. A basic result is then  the following:

\begin{prop}
\label{Siegelorbitboundprop}
For any Hodge-theoretic Siegel set $\mathfrak{G}$ associated to a Hodge structure $h$ on a polarized integral lattice $L$, a $\mathbb{Q}$-algebraic subvariety $X$ of the affine space associated to $L$, and orbit $\mathfrak{G} \cdot v$ of a Hodge vector $v$ for $h$, there exists a constant $c > 0$ such that
\[ \#\left( (\mathfrak{G} \cdot v) \cap X(\mathbb{Q}) \cap L \frac{1}{\nu} \right) \lesssim c \cdot \nu^{\min\{m + s - 2,\, m + r - 2,\, 2 \dim X\}} , \]
where $(r,s)$ is the signature of the underlying polarized lattice. 
\end{prop}
We note that the above result uses both that the Siegel set is Hodge-theoretic and that $v$ is a Hodge vector in a crucial way, and we expect that in general the intersections on the left are infinite without these assumptions. As a corollary we obtain:
\begin{cor}
\label{NLlocusbound}
For any polarized integral variation of Hodge structure $\mathbb{V}$ and orbit $O$ of $\mathbf{G}_{S}$, there exists a constant $c$ depending on $(O, \mathbb{V}, Q)$ such that
\[ \deg \NL_{O,\nu} \lesssim c \cdot \nu^{\min\{m + r - 2,\, m + s - 2,\, 2 \dim O\}} \]
where $(r,s)$ is the signature of $(\mathbb{V}, Q)$. 
\end{cor}

\noindent In the above result $\NL_{O,\nu}$ is the Noether-Lefschetz locus associated to $O$ and $\nu$, defined as 
\[ \NL_{O,\nu} = \left\{ s \in S(\mathbb{C}) :  \begin{array}{c}\textrm{ there exists a Hodge vector }v \in \mathbb{V}_{s} \\ \textrm{ with }v \in O(\mathbb{Q}), \, \, \nu(v) | \nu \end{array} \right\} . \]
Finally, we also have an overall bound on the Noether-Lefschetz locus for an arbitrary polarized VHS: 
\begin{cor}[\textbf{General Upper Bound}]
\label{NLqbound}
For any polarized integral variation of Hodge structure $\mathbb{V}$ there exists a constant $c$ depending only on $(\mathbb{V}, Q)$ such that
\[ \deg \NL_{q} \lesssim c \cdot q^{m^2 + \textrm{min}\{ r^2 + s^2, 2rs \} - 2} , \]
where $(r,s)$ is the signature of $(\mathbb{V}, Q)$. 
\end{cor}
\noindent In the above we define $\NL_{q}$ by
\[ \NL_{q} := \left\{ s \in S(\mathbb{C}) : \textrm{ there exists a Hodge vector }v \in \mathbb{V}_{s} \textrm{ with }\overline{Q}(v) = q \right\} . \]
The bound in \autoref{NLqbound} is deduced from \autoref{NLlocusbound} via a trick and we expect it can be substantially improved.

In explicit situations all of the above abstract results can be sharpened, and this is what we have done in the examples in the previous section. What happens then is that one knows explicitly the groups $\mathbf{G}_{S}$, $\mathbf{H}_{S}$, and $\mathbf{M}$, the structure of the Siegel sets becomes more explicit, and correspondingly the point counting becomes easier. We discuss this further below.

\subsection{A Sketch of the Ideas behind the Arguments}

\subsubsection{Reduction to Point Counting}
\label{redtoptcountsec}

Let us denote by $\pi : D_{S} \to \Gamma \backslash D_{S}$ the natural projection, and fix a Hodge subdatum $(\mathbf{M}, D_{M}) \subset (\mathbf{G}_{S}, D_{S})$. Suppose one wanted to merely count loci in $\Gamma \backslash D_{S}$ which are ``generalized Hecke translates'' of the image of some $D_{M} \subset D$: i.e., images under $\pi$ of $D_{M'}$, where $(\mathbf{M}', D_{M'})$ is a Hodge subdatum of $(\mathbf{G}_{S}, D_{S})$ which is $\mathbb{Q}$-coset equivalent to $(\mathbf{M}, D_{M})$. Then as we have discussed, one can view each such generalized Hecke translate as corresponding to a rational point of $\mathbf{G}_{S}/\mathbf{M}$. One can then try to pick a fundamental set for the action of $\Gamma$ on $\mathbf{G}_{S}(\mathbb{R})$, and count generalized Hecke translates by counting rational points in $\Gamma \backslash G(\mathbb{R})/\mathbf{M}(\mathbb{R})$. (The space $\Gamma \backslash G(\mathbb{R})/\mathbf{M}(\mathbb{R})$ does not usually have a reasonable topology, so one has to work directly with fundamental sets, but we ignore this for the moment.)

However in general one merely has a map $S \to \Gamma \backslash D_{S}$, and to make the strategy work one faces at least two potential problems.
\begin{itemize}
\item[(1)] The loci $\pi(D_{M'})$ might not intersect $\varphi(S)$. More generally, if one looks at the locus of $g \mathbf{M}(\mathbb{R}) \in \mathbf{G}_{S}(\mathbb{R})/\mathbf{M}(\mathbb{R})$ for which $\pi(g D_{M})$ intersects $\varphi(S)$ this locus might have empty interior, making it difficult to produce rational points inside it.
\item[(2)] It could be that many different varieties of the form $\pi(g D_{M})$, in particular infinitely many, intersect $\varphi(S)$ in the same locus. In other words, one has many different rational points of $\mathbf{G}_{S}/\mathbf{M}$ corresponding to the same Hodge locus, leading to overcounting.
\end{itemize}

The first problem (1) is solved using by the ``$\mathbb{V}$-likely hypothesis'', which can be used to guarantee intersections, following \cite{arXiv:2303.16179}. Then it would in fact follow from the Hodge-theoretic Zilber-Pink conjecture that (2) does not occur (more precisely, is controlled by a finite action of the normalizer of $\mathbf{M}$ in $\mathbf{G}_{S}$) away from a strict Zariski closed locus in $S$. This being unavailable, we need to assume the ``strongly $\mathbb{V}$-likely'' hypothesis, which, following the arguments in \cite{arXiv:2303.16179}, suffices to demonstrate the predictions of Hodge-theoretic Zilber-Pink in this case. 

To go from merely ``counting Hodge loci'' to a precise degree estimate one replaces $S$ with an appropriately chosen hyperplane section so that the Hodge loci of interest are zero dimensional. One then carefully partitions $S$ into analytic neighbourhoods $B \subset S(\mathbb{C})$ such that a Hodge locus in $B$ of the form $B \cap \varphi^{-1}(g \cdot D_{M})$ always consists of a single point. That this can be done for $B$ belonging to a finite definable analytic cover $S(\mathbb{C}) = \bigcup_{i=1}^{n} B_{i}$ is a consequence of the definability in $\mathbb{R}_{\textrm{an,exp}}$ of Hodge-theoretic period maps proven in \cite{defpermap}. Then by counting Hodge loci in each $B_{i}$ and summing over all $i$, one recovers an estimate for the number of points in the Hodge locus, which, having reduced to the zero dimensional case, is the degree estimate we wanted.

The task thus reduces to counting Hodge loci in each $B_{i}$. Following the theory in \cite{defpermap}, we can reduce to the case where each $B_{i}$ maps into a Hodge-theoretic Siegel set. After this, one establishes a correspondence between the Hodge locus points of interest and vectors in a Siegel set orbit, and counts points in the manner discussed above.

\subsubsection{Counting Points}

For the actual point counting itself, there are two tools. The first is the work of Gorodnik-Oh \cite{zbMATH05902967}, which gives general techniques based on adelic equidistribution results for counting points in subsets of homogeneous spaces for algebraic groups. These estimates are often sharp, and describe the asymptotic point counts using Euler-product-like expressions over all places of $\mathbb{Q}$. One can essentially always apply such techniques to count rational points inside compact subsets of $O(\mathbb{R})$, which leads to the lower bounds in \autoref{NLlocidegthmlower} and \autoref{abidemthm}. However it is unclear how to apply such techniques to count rational points in a general Siegel set orbit, especially given that usually such orbits are non-compact and have infinite volume in the natural $\mathbf{G}_{S}(\mathbb{R})$-invariant measure.

The upper bounds, therefore, are computed via a different method. We study carefully the geometry of Siegel set orbits, and reduce the problem to an analysis of orbits of certain ``standard Siegel sets''. Then we compute the orbits in explicit coordinates, and prove that one can bound the size of rational points in such orbits by the denominators of specially chosen coordinate entries. This in particular implies that rational points in Siegel set orbits whose denominators divide $\nu$ lie inside a compact subset of volume polynomially bounded by $\nu$. After this, the result follows from elementary bounds on the number of lattice points in a compact region. 

\subsection{Comparison with Other Work}
\label{compwithotherworksec}

Our results are similar in spirit to recent work studying the equidistribution of Hodge loci, in particular the papers \cite{zbMATH07201151} and \cite{zbMATH07643474}. There are a few differences and similarities.

\vspace{0.5em}

The first is that we estimate the degree of Hodge loci (or Noether-Lefschetz loci) for a polarized variation of Hodge structure $(S, \mathbb{V})$, while \cite{zbMATH07201151} and \cite{zbMATH07643474} both essentially estimate the degree of a corresponding ``locus of Hodge classes'' of $(S, \mathbb{V})$. Here we are borrowing the terminology of \cite{CDK} and Voisin \cite{voisin2010hodge} to describe the analytic loci in the vector bundle $\mathcal{H}^{\textrm{an}} = \mathbb{V} \otimes_{\mathbb{Z}} \mathcal{O}_{S^{\textrm{an}}}$ consisting of points $(s, v)$ with $v$ a Hodge class for the Hodge structure on $\mathbb{V}_{s}$. The vector bundle $\mathcal{H}^{\textrm{an}}$ has a canonical algebraic model $\mathcal{H}$ for which the components of the locus of Hodge classes are algebraic \cite{CDK}. That \cite{zbMATH07643474} count a (subset of the) Hodge class locus rather than the Noether-Lefschetz locus in $S$ can be seen in the statements of the asymptotic theorems \cite[Thm. 1.6, Thm. 1.7, Thm. 1.21, Thm. 1.22]{zbMATH07643474}. The situation in \cite[Thm. 1.1]{zbMATH07201151} is similar, where points in the Hodge locus are ``counted with multiplicity'', which in practice amounts to counting points by the number of Hodge classes which define them (but regarding Hodge classes differing by a scalar as equivalent).

Because a Noether-Lefschetz locus is always a projection under $\mathcal{H} \to S$ of a locus of Hodge classes, asymptotic lower bounds for degrees of Noether-Lefschetz loci imply asymptotic lower bounds for the degrees of Hodge class loci, and asymptotic upper bounds for the degrees of Hodge class loci imply asymptotic upper bounds for the degrees of Noether-Lefschetz loci. But there is in general a gap between the two, and understanding the size of this gap and when it occurs is the subject of non-trivial questions in unlikely intersection theory. For instance, the Zilber-Pink conjecture of \cite{BKU} predicts that if one considers a Hodge-generic family $f : X \to S$ of smooth projective degree $d \geq 5$ hypersurfaces in $\mathbb{P}^3$ with $S$ a curve, the Noether-Lefschetz locus of $S$ is a finite subset of $S(\mathbb{C})$. Formulated for the locus of Hodge classes, or when counting intersections with Noether-Lefschetz components in $U_{d}$ with multiplicity, the statement is false: if $s_{0} \in S$ is a point mapping to the Fermat surface in $U_{d}$ (recall \S\ref{NLexsec}), the locus of Hodge classes above $s_{0}$ consists of all points in $\mathbb{V}_{\mathbb{Z},s_{0}}$ belonging to an infinite lattice of large rank (c.f. \cite[Thm. 1]{zbMATH07074707}), and counting with multiplicity each $1$-dimensional subspace is counted separately. And situations where the Hodge locus of $S$ is infinite and the two types of loci should have different asymptotics also occur: if one takes $S$ to be a Shimura curve in $\mathcal{A}_g$ with $g \geq 2$, every Hodge locus point $s_0 \in S$ arises in infinitely many different ways as the intersection of $S$ with a special locus in $\mathcal{A}_g$, and from this one can show that counting the Hodge class locus (or counting with the multiplicity of the intersections) gives an overcount of the Hodge locus of $S$. (This phenomenon is also what we identified as ``potential problem (2)'' in \S\ref{redtoptcountsec}.)

\vspace{0.5em}

The second difference is that our upper bounds are global in nature. Although one could in principle use the methods of \cite{zbMATH07643474} to give upper-bound asymptotics for Hodge loci inside small neighbourhoods of $S(\mathbb{C})$ of uncontrolled size, to deduce a global point count on all of $S(\mathbb{C})$ one must, at least in the case of non-compact $S(\mathbb{C})$, sum over infinitely many such neighbourhoods. Without uniformly controlling the rates of convergence on all neighbourhoods at once, it is then difficult to obtain a global estimate. Our methods allow one to reduce the problem to counting rational points associated to finitely many Siegel sets, which crucially uses the limit theory of variations of Hodge structures. This finiteness ensures one can safely combine the asymptotic estimates in each region. An exception to this is \cite[Thm. 1.1]{zbMATH07201151}, which also achieves a global upper bound for variations of K3 type. 

\vspace{0.5em}

On the other hand, the underlying principle in both our approaches is the same. The idea is to use Hodge theory to replace degree estimates with counting problems for rational vectors. Over a compact set, explicit estimates may then be deduced from the work of Eskin, Gorodnik, Oh and others on such counting problems. We also both rely on measure theory, though we do our computations in coset spaces $(\mathbf{G}_{S}/\mathbf{M})(\mathbb{R})$ whereas \cite{zbMATH07643474} instead integrates a differential form on $S$ obtained by applying a ``push-pull'' procedure to a form on such a coset space. Ultimately such volume computations should amount to the same thing and we believe that one could use the results of \cite{zbMATH07643474} to deduce the analogues of the lower bound statements \autoref{NLlocidegthmlower} and \autoref{abidemthm} for an associated Hodge class locus (indeed, for \autoref{NLlocidegthmlower} this is basically contained in \cite[Thm. 1.6, Thm. 1.7]{zbMATH07643474}, and for split Jacobians there is \cite[Thm 1.14]{zbMATH07643474}). But for producing a lower bound for the Noether-Lefschetz locus itself, we do not know how to do so without additionally using arguments like those in \cite{arXiv:2303.16179} to solve the relevant unlikely intersection problems, and this is the approach to Noether-Lefschetz lower bounds we take in this paper. 

\subsection{Acknowledgements}

The author thanks both Salim Tayou and Nicolas Tholozan for several discussions about their work, and for explaining to him their perspective on Hodge locus asymptotics. He thanks Salim Tayou in particular for discussions at the MSRI (now SLMath) and at the third JNT Biennial Conference in Cetraro, and for pointing him to his paper \cite{zbMATH07201151}.

\section{Recollections on Period Maps}

\subsection{Siegel Sets}
\label{siegesetsec}

Let $(G, D) = (\mathbf{G}_{S}, D_{S})$ be the generic Hodge datum associated to $\mathbb{V}$, and $\varphi : S \to \Gamma \backslash D$ the period map induced by $\mathbb{V}$ with $\Gamma \subset G(\mathbb{Q})$ a neat arithmetic lattice containing the image of monodromy. By embedding $D$ as an open subvariety of its compact dual $\ch{D}$ the set $D$ inherits a natural $\mathbb{R}_{\textrm{alg}}$-definable structure. 

To define a Siegel set of $D$ we require some setup. Each such set will be defined as a certain orbit $\mathfrak{O} \subset D$ associated to a maximal compact subgroup $K \subset G(\mathbb{R})$ and a minimal parabolic $\mathbb{Q}$-subgroup $P \subset G$. The maximal compact group will be defined using a Hodge structure $o \in D$ by the following lemma:

\begin{lem}
\label{uniquemaxcompactlem}
Suppose that $o \in D$ is a point, and let $M_{o} \subset G(\mathbb{R})$ be its stabilizer. Then there is a unique maximal compact subgroup $K_{o} \subset G(\mathbb{R})$ containing $M_{o}$. 
\end{lem}

\begin{proof}
The Hodge structure $o : \mathbb{U} \to G_{\mathbb{R}}$ induces a grading $\mathfrak{g}_{\mathbb{C}} = \bigoplus_{i} \mathfrak{g}^{i,-i}$ through the adjoint action, and we define $K_{o}$ to be the exponential of the real Lie subalgebra $\mathfrak{k}_{o}$ underlying $\mathfrak{k}_{o,\mathbb{C}} = \bigoplus_{i\textrm{ even}} \mathfrak{g}^{i,-i}$. As the Lie algebra $\mathfrak{m}_{o}$ of $M_{o}$ may be identified with the $\mathfrak{g}^{0,0}_{\mathbb{R}}$ summand one clearly has $M_{o} \subset K_{o}$. A polarization on the Hodge structure $o$ induces a polarization on each simple summand $\mathfrak{u}$ of $\mathfrak{g}$ which is $\textrm{ad}\, \mathfrak{u}$-invariant when restricted to that summand, see \cite[III.A.5,  pg.75]{GGK}. Necessarily such a form is proportional to the Killing form $B$ on each summand, which implies that the restriction of the Killing form to $\mathfrak{k}_{o}$ is definite of a common sign, and hence negative definite because the Killing form is negative definite on $\mathfrak{m}_{o}$. This implies the restriction of the Killing form must be definite with a common (opposite) sign on the odd part of the Lie algebra. It then follows from \cite[\S4]{zbMATH03177763} that $\mathfrak{k}_{o}$ is a maximal compact subgroup containing $\mathfrak{m}_{o}$. That it is unique, which amounts to the statement that any other subalgebra of $\mathfrak{g}$ containing $\mathfrak{m}_{o}$ on which $B$ is negative definite lies in $\mathfrak{k}_{o}$, follows from the Hodge-Riemann bilinear relations for the polarizing form on $\mathfrak{g}$.
\end{proof}

\begin{lem}
\label{Sexistlem}
For any minimal parabolic $\mathbb{Q}$-subgroup $P$ of a reductive $\mathbb{Q}$-group $G$ and maximal compact subgroup $K \subset G(\mathbb{R})$, there exists a unique real torus $S_{P,K} \subset P_{\mathbb{R}}$ satisfying:
\begin{itemize}
\item[(i)] the torus $S_{P,K}$ is $P(\mathbb{R})$-conjugate to a maximal $\mathbb{Q}$-split torus of $P$; and
\item[(ii)] $S_{P,K}$ is stabilized by the Cartan involution associated to $K$. 
\end{itemize}
\end{lem}

\begin{proof}
This is \cite[Lem 2.1]{zbMATH06880895}. 
\end{proof}

\begin{defn}
\label{Atdef}
Suppose that $P \subset G$ is a minimal parabolic $\mathbb{Q}$-subgroup, and let $S = S_{P,o}$ be as in \autoref{Sexistlem}. For any real number $t > 0$ write 
\[ A_{t} = \{ \alpha \in S(\mathbb{R})^{+} : \chi(\alpha) \geq t \textrm{ for all }\chi \in \Delta \} , \]
where $\Delta$ is the set of simple roots of $G$ with respect to $S$, using the ordering induced by $P$. Then we define a Siegel set $\mathfrak{S} \subset G(\mathbb{R})$ (associated to $(P,o,t)$) to be a set $\Omega A_{t} K_{o}$ where $\Omega \subset P(\mathbb{R})$ is compact. We likewise define a Siegel set $\mathfrak{O} \subset D$ (associated to $(P,o,t)$) to be an orbit $\Omega A_{t} K_{o} \cdot o$.
\end{defn}

\subsection{Definable Period Maps}

We now recall one the main results of \cite{defpermap} (in the form corrected by \cite{zbMATH07720398}), which says roughly that local period maps associated to polarized integral variations of Hodge structures land inside Siegel sets. We also explain a very mild strengthening of this result which we will use in our arguments. 

We note that if $\mathfrak{F}$ is any definable fundamental set for $\Gamma$, there is a unique induced definable structure on $\Gamma \backslash D$ for which the map $\mathfrak{F} \to \Gamma \backslash D$ is definable; see \cite[Prop. 2.3]{zbMATH07855543}. The following is a minor variant of \cite[Thm. 1.5]{defpermap}, whose argument we summarize for completeness.

\begin{prop}[\cite{defpermap} + $\ep$]
\label{landsinsiegelset}
There exists a definable fundamental set $\mathfrak{F} \subset D$ for $\Gamma$, obtained as a finite union of Siegel sets $\mathfrak{O}_{i}$, such that the map $\varphi : S \to \Gamma \backslash D$ is $\mathbb{R}_{\textrm{an,exp}}$-definable for the induced definable structure on $\Gamma \backslash D$. In particular, there exists a finite cover $S = \bigcup_{i = 1}^{n} B_{i}$ by definable simply-connected opens and $\mathbb{R}_{\textrm{an,exp}}$-definable local lifts $\psi_{i} : B_{i} \to \mathfrak{O}_{i} \subset \mathfrak{F}$ of $\varphi$.

Fix a set of points $\mathcal{S} \subset S(\mathbb{C})$. Then we may choose the above data such that for each $B_{i}$ intersecting the topological closure $\mathcal{S}^{\textrm{top}}$ of $\mathcal{S}$ in $S(\mathbb{C})$, the Siegel set $\mathfrak{O}_{i} = \Omega_{i} A_{t,i} K_{o_{i}} \cdot o_{i}$ is associated to a point $o_{i} \in \psi_{i}(B_{i} \cap \mathcal{S})$. 
\end{prop}

\begin{proof}
We start by summarizing the proof in \cite{defpermap} of the first paragraph in the statement of \autoref{landsinsiegelset}. We then explain how to modify the argument to involve $\mathcal{S}$ at the end.

Applying Hironaka's Theorem we may reduce to the case where $S$ is the complement of a simple normal crossing divisor $E$ in a smooth projective variety $\overline{S}$. By passing to a finite \'etale cover we may assume monodromy at infinity is unipotent.  

It is clear we can construct such a $\psi$ in a small enough neighbourhood of any $s \in S$, so it suffices to show this on a small enough punctured neighbourhood $B = \Delta^{a} \times (\Delta^{\circ})^{b}$ of the simple normal crossing divisor $E$; this being done, the finiteness of the cover will follow from compactness. To do this we follow \cite[\S4]{defpermap}. Without loss of generality we take $a = 0$ by allowing factors with trivial monodromy. We let $\exp : \mathbb{H}^{b} \to (\Delta^{\circ})^{b}$ be the standard universal covering and let $\mathfrak{H} \subset \mathbb{H}$ be the standard Siegel fundamental set as in \cite[\S4.2]{defpermap}. Then we may construct a map $\widetilde{\psi} : \mathfrak{H}^{b} \to D$ by lifting $\varphi$, and define $\psi$ by inverting $\restr{\exp}{\mathfrak{H}^{b}}$ after choosing a branch cut and and composing with $\widetilde{\psi}$. By expanding the map $\restr{\exp}{\mathfrak{H}^{b}}$ in terms of its real-analytic components one sees it is $\mathbb{R}_{\textrm{an,exp}}$-definable, so to show that $\psi$ is $\mathbb{R}_{\textrm{an,exp}}$-definable it suffices to check the definability of $\widetilde{\psi}$.

To do this we apply the Nilpotent orbit theorem \cite{Schmid1973}, which tells us that 
\[ \widetilde{\psi} = \exp\left(\sum_{i=1}^{b} z_{i} N_{i} \right) \cdot \psi_{\textrm{lim}}(\exp(z)) \]
where $z = (z_{1}, \hdots, z_{b})$ are the natural coordinates for $\mathfrak{H}^{b}$, the $N_{i}$ are nilpotent operators, and $\psi_{\textrm{lim}}$ is an analytic function on $\Delta^{b}$, hence definable. Then as the left factor is just a polynomial in $z$, the result follows.

Following \cite[\S4.5]{defpermap}, it now suffices to check that the image of $\widetilde{\psi}$ lies in a finite union of Siegel sets. The monodromy representation associated to $\mathbb{V}$ gives a faithful representation $\rho : G \to \SL(V)$ where $V$ is an integral lattice, and hence induces a map $\iota : D \to D_{V}$, where $D_{V}$ is the symmetric space of positive-definite symmetric forms on $V_{\mathbb{R}}$; here the the map $\iota$ is given by sending the polarized Hodge structure $o$ to the positive-definite symmetric form $Q_{o} : (u,v) \mapsto Q_{\mathbb{C}}(C_{o} u, v)$, where $C_{o}$ is the Weil operator associated to $o$. From \cite[Prop 3.4]{zbMATH07778538} one learns that the inverse image under $\iota$ of a Siegel set of $D_{V}$ associated to $\widetilde{K}_{o}$ is contained in finitely many Siegel sets of $D$ for $K_{o}$. Thus we may reduce to the same claim for the image of $\iota \circ \widetilde{\psi}$. 

Using the explicit description of Siegel sets for $D_{V}$ in terms of the Gram-Schmidt process, the condition that the set $\mathfrak{B} = \textrm{im}(\iota \circ \widetilde{\psi}) \subset D_{V}$ lie inside a Siegel for $D_{V}$ is equivalent to the following claim: there exists a basis $\mathcal{E} = \{ e_{1}, \hdots, e_{m} \}$ for $V$ and a constant $\kappa$ such that the following inequalities hold for all $b \in \mathfrak{B}$:
\begin{itemize}
\item[(i)] $|b(e_{i}, e_{j})| < \kappa \,b(e_{i}, e_{i})$ for all $i, j$;
\item[(ii)] $b(e_{i}, e_{i}) < \kappa\, b(e_{j}, e_{j})$ for $i < j$; and
\item[(iii)] $\prod_{i} b(e_{i}, e_{i}) < \kappa \det b$.
\end{itemize}
The corresponding Siegel subset of $D_{V}$ may be taken to be associated to the maximal compact group $\widetilde{K}_{o}$ associated to a chosen point $o \in D_{S}$ (c.f. \cite[Thm. 3.3]{zbMATH07778538} and its proof). In particular we can take $o \in \psi(B)$.

We now check conditions (i), (ii) and (iii) by summarizing the argument in \cite[\S4.5]{defpermap}, which we refer to for details. Using the limiting weight filtration one can choose a basis $e_{1}, \hdots, e_{m}$ so that condition (iii) holds for some $\kappa$ as a consequence of \cite[Thm. 4.8]{defpermap}, and (ii) can be assumed for the same $\kappa$ by reordering the basis and partitioning the image as necessary. Schmid's $1$-dimensional result shows all of these conditions after restricting to any curve $\tau$ in $\mathfrak{H}^{b}$ for a constant $\kappa_{\tau}$ depending only on $\tau$, and so in particular shows (1) on such curves. One can then make the constant $\kappa_{\tau}$ appearing in (1) independent of $\tau$ using the fact that the coordinates of $Q_{\mathbb{C}}(C_{\widetilde{\psi}(z)}(-), (-))$ are ``roughly polynomial'' (see \cite[Lem. 4.5]{defpermap}) in the coordinates $z$ as a consequence of the results \cite{zbMATH04000092} and \cite{zbMATH03955093} on the asymptotics of Hodge forms. This completes the proof of the first paragraph of \autoref{landsinsiegelset}.

\vspace{0.5em}

Now we explain how to involve the set $\mathcal{S}$. The above argument has showed how, around any point $s \in \overline{S}(\mathbb{C})$, and given any ball $\overline{B} \subset \overline{S}(\mathbb{C})$ centred at $s$ with $\overline{B} = \Delta^{a + b}$, we can construct a period map $\psi$ on $B = \overline{B} \cap S(\mathbb{C})$ landing inside a Siegel set $\mathfrak{G} = \mathfrak{G}(s)$ depending on $s$. (Up to shrinking $B$ and taking a branch cut.) This in particular applies to points $s \in \overline{\mathcal{S}}^{\textrm{top}}$, where we write $\overline{\mathcal{S}}^{\textrm{top}}$ for the topological closure of $\mathcal{S}$ in $\overline{S}(\mathbb{C})$. Now as discussed above, the maximal compact of the resulting Siegel set can be chosen to be associated to any point $o \in \psi(B)$. Since $\mathcal{S}$ is dense in $\overline{\mathcal{S}}^{\textrm{top}}$, we can in particular choose $o \in \psi(B \cap \mathcal{S})$.

The set $\overline{\mathcal{S}}^{\textrm{top}}$ is a closed subset of a compact set, so is in particular compact. Thus, after applying this argument at every point $s \in \overline{\mathcal{S}}^{\textrm{top}}$ and obtaining an appropriate cover $\{ \overline{B}_{s} \}_{s \in \overline{\mathcal{S}}^{\textrm{top}}}$ of $\overline{\mathcal{S}}^{\textrm{top}}$, we can find a finite subcover. This finite subcover may then be extended to a finite cover of all of $\overline{S}$ which induces the cover of $S$ given in the statement. 
\end{proof}

\section{Siegel sets and their Orbits}
\label{Siegelsetsec}

In this section we study carefully the structure of, and point counting in, Siegel set orbits for the group $\SO(r,s)$. 

\subsection{Standard Siegel Sets}
\label{stdsiegelset}

We set $m = r+s$, and $G = \SO(r,s)$. We will view $G$ as the stabilizer of the bilinear form $Q = -I_{r} \oplus I_{s}$. Replacing $Q$ with $-Q$ and reordering the standard coordiantes $x_{1}, \hdots, x_{m}$ if necessary, we may assume that $r \geq s$; none of our results on Siegel set orbits will depend on this convention. 

We will assume that $\mathfrak{S}$ is a Siegel set constructed with respect to a certain special choice of maximal compact subgroup $K$ and parabolic subgroup $P$; we justify this in \autoref{Siegelconjlem} below by showing that one can always reduce to this setting by a change of coordinates. We fix the maximal compact subgroup $K = S(O(r) \times O(s))$. To construct $P$ we introduce a second set of coordinates, $w_{1}, \hdots, w_{m}$, which are related to the coordinates $x_{1}, \hdots, x_{m}$ by 
\begin{align}
w_{i} &= x_{i} - x_{r+i} & 1 \leq i \leq s \\
w_{i} &= x_{i} & s+1 \leq i \leq r \\
w_{i} &= x_{i} + x_{i-r} & r+1 \leq i \leq m .
\end{align}
The change of basis matrices from $\overline{x}$-coordinates to $\overline{w}$-coordinates are given by 
\[ S = \begin{pmatrix} I_{s} & 0 & -I_{s} \\ 0 & I_{r-s} & 0 \\ I_{s} & 0 & I_{s} \end{pmatrix} , \hspace{2em} S^{-1} = \frac{1}{2} \begin{pmatrix} I_{s} & 0 & I_{s} \\ 0 & 2 I_{r-s} & 0 \\ -I_{s} & 0 & I_{s} \end{pmatrix} . \]

And in the new coordinate frame, the bilinear form $Q$ is given by
\begin{align*}
Q_{S} &:= \frac{1}{4} \begin{pmatrix} I_{s} & 0 & -I_{s} \\ 0 & 2 I_{r-s} & 0 \\ I_{s} & 0 & I_{s} \end{pmatrix} \begin{pmatrix} - I_{s} & 0 & 0 \\ 0 & -I_{r-s} & 0 \\ 0 & 0 & I_{s} \end{pmatrix} \begin{pmatrix} I_{s} & 0 & I_{s} \\ 0 & 2 I_{r-s} & 0 \\ -I_{s} & 0 & I_{s} \end{pmatrix} \\
&= \frac{1}{4} \begin{pmatrix} I_{s} & 0 & -I_{s} \\ 0 & 2 I_{r-s} & 0 \\ I_{s} & 0 & I_{s} \end{pmatrix} \begin{pmatrix} - I_{s} & 0 & -I_{s} \\ 0 & -2 I_{r-s} & 0 \\ -I_{s} & 0 & I_{s} \end{pmatrix} \\
&= \frac{1}{4} \begin{pmatrix} 0 & 0 & -2 I_{s} \\ 0 & -4 I_{r-s} & 0 \\ -2 I_{s} & 0 & 0 \end{pmatrix} .
\end{align*}
We compute the Lie algebra in the $\overline{w}$-coordinate system, which is given by matrices $W$ satisfying $W^{t} Q_{S} + Q_{S} W = 0$. For such $W$, one calculates that 
\begin{align*}
0 &= \begin{pmatrix} W^{t}_{11} & W^{t}_{21} & W^{t}_{31} \\ W^{t}_{12} & W^{t}_{22} & W^{t}_{32} \\ W^{t}_{13} & W^{t}_{23} & W^{t}_{33} \end{pmatrix} \begin{pmatrix} 0 & 0 & I_{s} \\ 0 & 2I_{r-s} & 0 \\ I_{s} & 0 & 0 \end{pmatrix} + \begin{pmatrix} 0 & 0 & I_{s} \\ 0 & 2I_{r-s} & 0 \\ I_{s} & 0 & 0 \end{pmatrix} \begin{pmatrix} W_{11} & W_{12} & W_{13} \\ W_{21} & W_{22} & W_{23} \\ W_{31} & W_{32} & W_{33} \end{pmatrix} \\
&= \begin{pmatrix} W^{t}_{31} + W_{31} & 2 W^{t}_{21} + W_{32} & W^{t}_{11} + W_{33} \\ W^{t}_{32} + 2 W_{21} & 2 (W^{t}_{22} + W_{22}) & W^{t}_{12} + 2 W_{23} \\ W^{t}_{33} + W_{11} & 2 W^{t}_{23} + W_{12} & W^{t}_{13} + W_{13} \end{pmatrix} .
\end{align*}
So in the $\overline{w}$-coordinate system one obtains the description
\begin{align*}
\mathfrak{so}(r,s)_{\overline{w}} = \left\{ \begin{pmatrix} W_{11} & - 2 W^{t}_{23} & W_{13} \\ W_{21} & W_{22} & W_{23} \\ W_{31} & - 2 W^{t}_{21} & -W^{t}_{11} \end{pmatrix} : W_{22} = -W^{t}_{22}, \hspace{0.5em} W_{31} = - W^{t}_{31}, \hspace{0.5em} W_{13} = - W^{t}_{13} \right\} .
\end{align*}

Write $U_{s} \subset \GL_{s}$ the group of $s \times s$ upper triangular matrices, and $\mathfrak{u}_{s}$ for its Lie algebra. We define a Lie algebra $\mathfrak{p}$ by the equations
\begin{align*}
\mathfrak{p} &= \{ W \in \mathfrak{so}(r,s)_{\overline{w}} : \hspace{0.5em} W_{21} = 0, \hspace{0.5em} W_{31} = 0, \hspace{0.5em} W_{11} \in \mathfrak{u}_{s}(\mathbb{R}) \} .
\end{align*}
In $\overline{w}$-coordinates, the corresponding matrices have the form
\begin{align*}
&= \left\{ \begin{pmatrix} Y & - 2 W^{t}_{23} & W_{13} \\ 0 & W_{22} & W_{23} \\ 0 & 0 & -Y^{t} \end{pmatrix} : Y \in U_{s}(\mathbb{R}), \hspace{0.5em} W_{13} = - W^{t}_{13} \right\} .
\end{align*}
We write $\mathfrak{n} \subset \mathfrak{p}$ for the nilpotent Lie algebra defined by setting $W_{22} = 0$ and mandating that the diagonal entries of $Y$ are all equal to $0$, and set $U = \textrm{exp}(\mathfrak{n})$. 

\begin{lem}
The set $\mathfrak{p}$ is the Lie algebra of a minimal parabolic subalgebra of $\mathfrak{g}$.
\end{lem}

\begin{proof}
\cite[\S2.3]{he2013certainunitarylanglandsvoganparameters} gives a description of a minimal parabolic subgroup of $G$. One immediately checks that our algebra $\mathfrak{p}$ and the group appearing in \cite[\S2.3]{he2013certainunitarylanglandsvoganparameters} have the same dimension, so it suffices to check that $\mathfrak{p}$ is parabolic. By \cite[Def. 1.1]{zbMATH05987928} it suffices to check that the orthogonal complement $\mathfrak{p}^{\perp}$ with respect to the Killing form is nilpotent, which is easily checked. 
\end{proof}

Write $P$ for the associated Parabolic subgroup of $G$. By exponentiating, we see that 
\begin{equation}
\label{Pmatform}
P \subset \left\{ \begin{pmatrix} A & B & C \\ 0 & D & E \\ 0 & 0 & A^{t,-1} \end{pmatrix} : A \in U_{s}(\mathbb{R}), B, C, E \in \GL_{m}(\mathbb{R}), D \in O(r-s)(\mathbb{R}) \right\} .
\end{equation}
We now describe the split torus $S \subset P$ determined by \autoref{Sexistlem}. The Cartan involution associated to $K$ in $\overline{x}$-coordinates is given by $\theta(X) = - X^{t}$. In the $\overline{w}$-coordinates this becomes
\begin{align*}
\theta(W) &= S \theta(S^{-1} W S) S^{-1} \\
&= - S S^{t} W^{t} S^{-1,t} S^{-1} .
\end{align*}
The matrices $S S^{t}$ and $S^{-1,t} S^{-1} = (S S^{t})^{-1}$ are both diagonal, so one sees that the diagonal matrices of $P$ are preserved under the involution. We then take the Lie algebra $\mathfrak{a}$ of $S$ to be the intersection in the $\overline{w}$-coordinates of $\mathfrak{p}$ and the diagonal subgelbra of $\mathfrak{gl}_{m}$. It is then clear that $S$ is split over $\mathbb{Q}$ of the correct dimension (compare with the torus in \cite[\S2.3]{he2013certainunitarylanglandsvoganparameters}), and stable under $\theta$. 

The rank of the group $\SO(r,s)$ is $s$ (recall we have assumed $r \geq s$), so we should have $s$ simple roots $\chi_{1}, \hdots, \chi_{s}$ matching \autoref{Atdef}. Since our torus $S$ lies in the diagonal subgroup $D_{m}$ of $\GL_{m}$, we may obtain these roots from restrictions of roots of $D$ acting on $\mathfrak{gl}_{m}$, in particular, we have 
\begin{align*}
A_{t} := \{ [a_{ij}] \in D_{m}(\mathbb{R})^{+} \cap \SO(r,s)(\mathbb{R}) : a_{ii} \geq t a_{(i+1)(i+1)} \hspace{0.5em} \textrm{ for } \hspace{0.5em} 1 \leq i \leq s-1, \hspace{1em} a_{s} \geq t \} .
\end{align*}

\begin{lem}
\label{Siegelconjlem}
Let $K' \subset G(\mathbb{R})$ be a maximal compact subgroup and $\mathfrak{S}'$ a Siegel set for $G(\mathbb{R})$ relative to $K'$ and a minimal paraoblic $P'$. Then there exists $\gamma \in G(\mathbb{Q})$, $\tau \in U(\mathbb{R})$ and $\beta \in S(\mathbb{R})$, and a Siegel set $\mathfrak{S}$ for $(P,S,K)$ such that $\gamma \mathfrak{S}' \gamma^{-1} \sigma \subset \mathfrak{S}$ where $\sigma = \tau \beta$. We have $\gamma P' \gamma^{-1} = P$ and $\sigma K \sigma^{-1} = \gamma^{-1} K' \gamma$.
\end{lem}

\begin{proof}
Our argument is a mild variant of \cite[Lem. 3.8]{zbMATH06880895}. Let $P'$ be the parabolic group associated to $\mathfrak{S}' = \Omega' A'_{t} K'$. Then as both $P$ and $P'$ are minimal parabolic $\mathbb{Q}$-subgroups, there exists $\gamma \in G(\mathbb{Q})$ such that $\gamma P' \gamma^{-1} = P$. Since $K$ and $\gamma K' \gamma^{-1}$ are both maximal compact subgroups, there exists $\sigma \in G(\mathbb{R})$ such that $\sigma K \sigma^{-1} = \gamma^{-1} K' \gamma$.

We now consider the Iwasawa decomposition of $G(\mathbb{R})$ with respect to the Cartan involution $\theta$. We recall that this is induced by a decomposition $\mathfrak{g} = \mathfrak{n} \oplus \mathfrak{a} \oplus \mathfrak{k}$, where $\mathfrak{q} = \mathfrak{n} \oplus \mathfrak{a}$ and $\mathfrak{k}$ are both eigenspaces for $\theta$, and $\mathfrak{a}$ is a choice of maximal abelian subalgebra of $\mathfrak{q}$. Since the Lie algebra of $S$ is stabilized by $\theta$, we may choose $\mathfrak{a}$ to be the Lie algebra of $S$. It is then directly checked that $\mathfrak{n}$ is a sum of root spaces for $S$. It follows that $G(\mathbb{R}) = U(\mathbb{R})S(\mathbb{R})^{+} K$, and that we may choose $\sigma = \tau \beta$ with $\beta \in S^{+}(\mathbb{R})$ and $\tau \in U(\mathbb{R})$. It follows that $(\gamma \sigma)^{-1} P' (\gamma \sigma) = P$ and hence by \autoref{Sexistlem} that $(\gamma \sigma)^{-1} S' (\gamma \sigma) = S$, so we may write $A_{t} = \sigma^{-1} \gamma^{-1} A'_{t} \gamma \sigma$. Finally, as in the proof of \cite[Lem. 3.8]{zbMATH06880895}
\begin{align*}
\gamma^{-1} \mathfrak{S}' \gamma \sigma &= (\gamma^{-1} \Omega' \gamma) \sigma \sigma^{-1} \gamma^{-1} A'_{t} \gamma \sigma \sigma^{-1} \gamma^{-1} K' \gamma \sigma \\
&= (\gamma^{-1} \Omega' \gamma) \tau \beta A_{t} K .
\end{align*}
so the result follows choosing $\Omega \supset (\gamma^{-1} \Omega' \gamma) \tau$ and $A_{s} \supset \beta A_{t}$ sufficiently large.
\end{proof}

\subsection{Orbits of Hodge Vectors}

To compare the above Siegel sets to the ones given in \S\ref{siegesetsec}, we fix a Hodge structure $h : \mathbb{U} \to \SO(r,s)(\mathbb{R})$, let $K_{h}$ be the associated maximal compact, and let $(P_{h}, S_{h}, K_{h})$ be an associated Siegel triple. As in \S\ref{siegesetsec}, the group $K_{h}$ is the intersection with $\SO(r,s)(\mathbb{R})$ of $O(Q_{h})$, where $Q_{h}(u,w) = Q(C_{h} u, w)$ is the Hodge-metric at $h$. 

\begin{notn}
Given a vector $v \in \mathbb{R}^{m}$, we write either $v_{i}$ or $[v]_{i}$ for its $i$'th entry. We let $\ell = \ell(v) \in \{ r+1, \hdots, m \}$ be the smallest index $i$ in the specified range for which $v_{i} \neq 0$. If we write $[v]_{\ell}$ we understand $[v]_{\ell(v)}$. If no such index exists, then we say $\ell = \ell(v)$ is undefined.
\end{notn}

\begin{prop}
\label{Hodgevectororbitbound}
Fix a Siegel set $\mathfrak{G}$ for $(P_{h}, S_{h}, K_{h})$ and a real Hodge vector $v$ for $h$. Then there exists a basis $\mathcal{F} = \{ f_{1}, \hdots, f_{m} \}$ for $\mathbb{Q}^{m}$ with the following properties:
\begin{itemize}
\item[(1)] in the basis $\mathcal{F}$, the form $Q$ is equal to $-I_{r} \oplus I_{s}$ and the special orthogonal group of $Q$ is equal to $\SO(r,s)$;
\item[(2)] writing $\overline{x}$ for the coordinates with respect to $\mathcal{F}$ and defining $\overline{w}$ as in \S\ref{stdsiegelset}, there exists real constants $c, \kappa > 0$, independent of $w \in \mathfrak{G} \cdot v$, such that
\begin{align}
\label{Siegeleq1}
| w_{\ell} |^{-1} + \kappa &> c |w_{i}| , & 1 \leq i \leq s , \\
\label{Siegeleq2}
\kappa &> |w_{i}| , & s+1 \leq i \leq m . 
\end{align}
if $\ell = \ell(w)$ is defined, and $\| w \| < \kappa$ otherwise. 
\end{itemize}
\end{prop}

\begin{proof}
Define $\gamma$ and $\sigma$ as in \autoref{Siegelconjlem} so that $\gamma \mathfrak{G} \gamma^{-1} \sigma \subset \mathfrak{S}$, where $\mathfrak{S}$ is one of the standard Siegel sets constructed above. We define $\mathcal{F}$ as $\gamma \cdot \mathcal{E}$, where $\mathcal{E} = \{ e_1, \hdots, e_m \}$ is the standard basis. Then we have 
\[ (\gamma \mathfrak{G} \gamma^{-1}) \gamma v \subset \mathfrak{S} \sigma^{-1} \gamma v . \]
The left-hand side is nothing more than the orbit $\mathfrak{G} \cdot v$ except in the coordinates determined by $\mathcal{F}$, so it suffices to prove the inequalities for the points inside the right-hand side. The vector $\sigma^{-1} \gamma v$ is a real Hodge vector for the Hodge structure $\sigma^{-1} \gamma \cdot h$. Moreover because $\sigma K \sigma^{-1} = \gamma^{-1} K_{h} \gamma$, the maximal compact associated to $\mathfrak{S}$ agrees with the stabilizer of the Hodge metric for $\sigma^{-1} \gamma \cdot h$. It thus suffices to prove the statement of the proposition when $\mathfrak{G} = \mathfrak{S}$ is a standard Siegel set which is also a Hodge-theoretic Siegel set for $h$.

We thus reduce to considering the orbit $\Omega A_{t} K \cdot v$, with $K$ as in the previous section. Let $\mathbb{C}^{m} = \bigoplus H^{p,q}$ be the Hodge decomposition associated to $h$, where we assume the Hodge structure has weight $w = p + q$ where $w = 2k$ is even. Define $V_{e} = \bigoplus_{p\textrm{ even}} H^{p,w-p}$ and $V_{o} = \bigoplus_{p\textrm{ odd}} H^{p,w-p}$, which we may view as subspaces of $\mathbb{R}^{m}$. The Weil operator $C_{h}$ acts as $i^{2p - w} = (-1)^{p-k}$ on $H^{p,q}$, so the forms $Q$ and the Hodge metric $Q_{h}$ agree on one of the two summands $\{ V_{e}, V_{o} \}$ and agree up to a sign on the other.

The maximal compact group associated to $h$ preserves both $V_{e}$ and $V_{o}$. Since we have reduced to the setting where the maximal compact associated to $h$ is $K = S(O(r) \times O(s))$, it follows that $\{ V_{e}, V_{o} \} = \{ \textrm{span}\{ e_{1}, \hdots, e_{r} \}, \textrm{span}\{ e_{r+1}, \hdots, e_{m} \} \}$. Then because $v$ is Hodge, hence $v \in V_{e}$, we observe that all the vectors in $K \cdot v$ lie inside $V_{e}$ as well. In particular, if $k v \in K \cdot v$, then $[k v]_{i} = \pm  [k v]_{r+i}$ for $1 \leq i \leq s$ in $\overline{w}$-coordinates, with the sign depending on which of the two summands in $\{ \textrm{span}\{ e_{1}, \hdots, e_{r} \}, \textrm{span}\{ e_{r+1}, \hdots, e_{m} \} \}$ the space $V_{e}$ is identified with. We will assume the sign is positive, with the other case handled analogously.

\vspace{0.5em}
We now consider a vector $w \in \Omega A_{t} K \cdot v$ which we may write as $w = \omega a k \cdot v$ with $k \in K$, $a \in A_{t}$, and $\omega \in \Omega$. After possibly adjusting $a$ and $\omega$ and decreasing $t$, we may assume that all diagonal entries of $\omega$ are equal to $1$ except for possibly those that lie in the central block. We write $a_{i}$ for the $i$'th diagonal entry of $a$ in the $\overline{w}$-coordinate system. We note that 
\[ [a k \cdot v]_{i} = a_{i} [k \cdot v]_{i} = a_{i} [k \cdot v]_{r+i} = a^2_{i} [a k \cdot v]_{r+i} \]
for $1 \leq i \leq s$. We may assume $\ell_{0} = \ell(k \cdot v)$ is defined: otherwise $[ak \cdot v] = [k \cdot v]$ and the entries of $[\omega a k \cdot v] = [\omega k \cdot v]$ are universally bounded as both $\omega$ and $k$ range over a compact set.

We then have
\begin{align*}
[a k \cdot v]^{-1}_{\ell_{0}} = [k \cdot v]^{-2}_{\ell_{0} - r} [a k \cdot v]_{\ell_{0}-r} .
\end{align*}
Choose $j$ such that $\ell_{0} \leq j \leq r+i$ and $[k \cdot v]_{j-r} \neq 0$. From the inequalities defining $A_{t}$ we then get
\begin{align*}
[a k \cdot v]^{-1}_{\ell_{0}} \geq t^{j-\ell_{0}} [k \cdot v]^{-2}_{\ell_{0} - r} [k \cdot v]^{-1}_{j-r} [a k \cdot v]_{j-r} ,
\end{align*}
We now take $c$ such that $0 < c \leq [k \cdot v]^{-1}_{\ell_{0} - r} [k \cdot v]^{-1}_{j - r} t^{j-\ell_{0}}$, which we can do independently of $w$ and $i$ since $[k \cdot v]^{-1}_{\ell_{0} - r} [k \cdot v]^{-1}_{j-r}$ is bounded from below independently of $w$. Taking $j = r+i$ in the above gives
\begin{align*}
[a k \cdot v]^{-1}_{\ell_{0}} \geq c [a k \cdot v]_{i} ,
\end{align*}
for all $1 \leq i \leq s$. One computes immediately from the description of $P$ in (\ref{Pmatform}) that $\ell(\omega a k \cdot v) = \ell(k \cdot v) = \ell_{0}$ and that $[\omega a k \cdot v]_{\ell_{0}} = [a k \cdot v]_{\ell_{0}}$; in particular $\ell(w) = \ell(\omega a k \cdot v)$ is defined. It thus follows that 
\begin{equation}
\label{boundbywell}
[w]^{-1}_{\ell} \geq c [a k \cdot v]_{i} ,
\end{equation}
for all $1 \leq i \leq s$.

\vspace{0.5em}

Now we also have $a_{s} \geq t$, and hence $a_{i} \geq t^{s-i+1}$ for all $1 \leq i \leq s$. This implies that $a_{r+i} = a^{-1}_{i} \leq t^{s-i+1}$ for all $1 \leq i \leq s$, in particular, there is an absolute bound $\kappa$ such that
\[ \kappa > [a k \cdot v]_{i} \]
for all $i \geq r+1$. We can even assume this is true for all $i \geq s+1$ using that $A_{t}$ acts trivially on the coordinates indexed by $s+1, \hdots, r$. 

\vspace{0.5em}

We conclude that for all $1 \leq i \leq m$ we have 
\begin{align*}
[w]^{-1}_{\ell} + \kappa &\geq c [a k \cdot v]_{i} , & 1 \leq i \leq s \\
\kappa &> [a k \cdot v]_{i} , & s+1 \leq i \leq m .
\end{align*}
Then since $\omega$ lies in a compact subset $\Omega \subset P(\mathbb{R})$, the inequalities continue to hold with $[\omega a k \cdot v]_{i}$ replacing $[a k \cdot v]_{i}$ after possibly adjusting $c$ and $\kappa$.
\end{proof}

\subsection{Point Counting Upper Bounds}
\label{pointcountingsec}

We continue with the notation and setup of the preceding section. We view $\mathbb{R}^{m}$ as the set of real points of the standard $m$-dimensional affine space $\mathbb{A}^{m}$.

\begin{prop}
\label{siegelorbitnubound}
Fix a Siegel set $\mathfrak{G}$ for $(P_{h}, S_{h}, K_{h})$ and a real Hodge vector $v$ for $h$. Then there exists a real constant $\rho$ such that
\[ \left| (\mathfrak{G} \cdot v) \cap \left( \frac{1}{\nu} \mathbb{Z}^{m} \right) \right| \lesssim \rho \cdot \nu^{m + s - 2} \]
as $\nu \to \infty$. 
\end{prop}

\begin{proof}
The statement is unchanged (up to adjusting $\rho$) after a $\mathbb{Q}$-linear change of coordinates, so we may work in the coordinate system $\overline{w}$ of \autoref{Hodgevectororbitbound}. In the region of $\mathfrak{G} \cdot v$ where $\ell(w)$ is undefined we then know that $\| w \|$ is absolutely bounded, and on this region the result follows by projecting to a coordinate hyperplane. We may thus assume the inequalities (\ref{Siegeleq1}) and (\ref{Siegeleq2}) hold for some fixed $\ell = \ell(w)$. Rearranging these inequalities one obtains
\begin{align}
\label{reSiegeleq1}
\kappa' &> c |w_{i}| |w_{\ell}| , & 1 \leq i \leq s , \\
\label{reSiegeleq2}
\kappa &> |w_{i}| , & s+1 \leq i \leq m ,
\end{align}
where we take $\kappa' > 1 + \kappa^2 > 1 + \kappa |w_{\ell}|$ using that $\kappa > |w_{\ell}|$ from (\ref{reSiegeleq2}). 

For each index $r+1 \leq \ell \leq m$ we define a map $\tau_{\ell} : \mathbb{A}^{m} \to \mathbb{A}^{m-1}$ via the coordinate functions 
\begin{equation}
\label{wprojcoords}
(w_{1} w_{\ell}, \hdots, w_{\ell-r-1} w_{\ell}, w_{\ell-r+1} w_{\ell}, \hdots, w_{s} w_{\ell}, w_{s+1}, \hdots, w_{m}) . 
\end{equation}
Letting $X_{\ell} \subset \mathbb{A}^{m}$ be the locus where $w_{\ell} \neq 0$, we observe that the induced map
\[ X_{\ell} \cap (\SO(r,s) \cdot v) \xrightarrow{\tau_{\ell}} \mathbb{A}^{m-1} \]
is injective. Indeed the locus $\SO(r,s) \cdot v$ is defined by an equation 

\begin{equation}
\label{SOorbiteq}
\sum_{i=1}^{s} w_{i} w_{r+i} + \sum_{j=1}^{r-s} w_{j}^2 = c
\end{equation}
for some constant $c$. Given the coordinates (\ref{wprojcoords}) and the condition that $w_{\ell} \neq 0$, one recovers by dividing by $w_{\ell}$ all the coordinates $w_{i}$ except $w_{\ell-r}$. But then one can determine each term appearing in (\ref{SOorbiteq}) except for $w_{\ell-r} w_{\ell}$, so one can again divide through by $w_{\ell}$ and recover $w_{\ell-r}$.

The problem reduces to bounding rational points in the image of $\tau_{\ell}$. From the inequalities (\ref{reSiegeleq1}) and (\ref{reSiegeleq2}) we see that $\tau_{\ell}(\mathfrak{G} \cdot v \cap X_{\ell}(\mathbb{R}))$ lies inside some compact set $\mathcal{K}$, and the rational points in $\frac{1}{\nu} \mathbb{Z}^{m}$ have image inside $\frac{1}{\nu^2} \mathbb{Z}^{s-1} \oplus \frac{1}{\nu} \mathbb{Z}^{m-s}$. The number of such points is then at most $\rho \cdot \nu^{2(s-1)} \nu^{m-s} = \rho \cdot \nu^{m + s-2}$ choosing $\rho$ large enough. 
\end{proof}

\begin{prop}
\label{genericsiegelorbitbound}
Fix a Siegel set $\mathfrak{G}$ for $(P_{h}, S_{h}, K_{h})$ and a real Hodge vector $v$ for $h$. Suppose that $Y \subset \mathbb{A}^{m}$ is an irreducible $\mathbb{Q}$-algebraic subvariety. Then there exists a real constant $c$, dependent on $Y$, such that
\[ \left| (\mathfrak{G} \cdot v) \cap Y(\mathbb{Q}) \cap \left( \frac{1}{\nu} \mathbb{Z}^{m} \right) \right| \lesssim c \cdot \nu^{2 \dim Y} \]
as $\nu \to \infty$. 
\end{prop}

\begin{proof}
Once again we work in the coordinate system $\overline{w}$ of \autoref{Hodgevectororbitbound}. We write $V_{i} \subset \mathbb{A}^{m}$ for the locus where $w_{i}$ vanishes, and let $X_{i} = \mathbb{A}^{m} \setminus V_{i}$ be the complement. Let $E \subset [m] = \{ 1, \hdots, m \}$ be the subset of indices $i$ such that $Y \subset V_{i}$. Then by arguing inductively on $Y \cap V_{j}$ for $j \in [m] \setminus E$, we may reduce to considering just those rational points in $Y' := Y \cap \bigcap_{j \in [m] \setminus E} X_{j}$. 

We may assume $([m] \setminus E) \cap \{ r+1, \hdots, m \}$ is non-empty. Otherwise $\ell(w)$ is undefined for each $w \in Y(\mathbb{R})$, and the desired result follows by projecting $(\mathfrak{G} \cdot v) \cap Y(\mathbb{Q})$ to a coordinate hyperplane and using that $\| w \|$ is absolutely bounded. Then we let $\ell$ be the smallest entry of $([m] \setminus E) \cap \{ r+1, \hdots, m \}$.

We consider a projection $\tau'_{\ell} : \mathbb{A}^{m} \to \mathbb{A}^{\dim Y}$ defined by composing the map $\tau_{\ell}$ defined in the proof of \autoref{siegelorbitnubound} with a further projection $\mathbb{A}^{m-1} \to \mathbb{A}^{\dim Y}$ onto a coordinate hyperplane and a translation $\mathbb{A}^{\dim Y} \to \mathbb{A}^{\dim Y}$ by an integral vector in $\mathbb{Z}^{\dim Y}$. We may choose these projections and translations such that $Y \to \mathbb{A}^{\dim Y}$ is quasi-finite away from a closed $\mathbb{Q}$-algebraic sublocus of $Y$ which we handle inductively. For the remaining points we are reduced by taking the image under $\tau'_{\ell}$ to counting points in $\frac{1}{\nu^2} \mathbb{Z}^{\dim Y}$ which lie in a compact region, so the result follows. 
\end{proof}

\begin{proof}[Proof of (\ref{Siegelorbitboundprop}):]
This follows from \autoref{siegelorbitnubound} and \autoref{genericsiegelorbitbound} after recalling that we have adopted the convention that $r \geq s$ throughout \S\ref{Siegelsetsec}, and replacing $Q$ by $-Q$ does not change the group preserving the form. 
\end{proof}

\section{Assorted Tools}

\subsection{A Definability Lemma}

We start with a basic lemma regarding definable sets, which follows from definable cell decomposition. A definable cell decomposition is a partition of $\mathbb{R}^{n}$ into definable subsets, called cells, defined inductively as follows (c.f. \cite[Ch. 3, 2.3]{zbMATH01160037}):
\begin{itemize}
\item[(i)] when $n = 1$, the cells are open intervals and singleton sets;
\item[(ii)] given a cell $D \subset \mathbb{R}^{n-1}$, a cell of $\mathbb{R}^{n}$ is either:
\begin{itemize}
\item[-] the graph of a continuous definable function $f : D \to \mathbb{R}$ viewed in the natural way as a subset of $\mathbb{R}^{n}$; or
\item[-] the region in $\mathbb{R}^{n}$ defined by
\[ \{ (x, y) \in D \times \mathbb{R} : f(x) < y < g(x) \} \] 
where $f, g : D \to \mathbb{R}$ are continuous definable functions satisfying $f < g$.
\end{itemize}
\end{itemize}
Then one has \cite[Ch. 3, 2.11]{zbMATH01160037}:
\begin{thm}
For any definable subset $F \subset \mathbb{R}^{n}$ there exists a cell decomposition of $\mathbb{R}^{n}$ such that $F$ is a finite union of cells. 
\end{thm}

We will find the following lemma useful in later arguments.

\begin{lem}
\label{finitepartlem}
Let $A$ and $B$ be definable sets, and $F \subset A \times B$ a definable subset such that the intersections $F \cap \{ a \} \times B$ are all finite. Then there exists a definable partition $F = F_{1} \sqcup \cdots \sqcup F_{k}$ such that the intersections $F_{i} \cap \{ a \} \times B$ have cardinality at most $1$ for each $a$ and $i$.
\end{lem}

\begin{proof}
We may reduce to the case where $A = \mathbb{R}^{n}$ and $B = \mathbb{R}^{m}$, so $F$ is a definable subset of $\mathbb{R}^{n+m}$. We then consider a cell decomposition of $(F, \mathbb{R}^{n+m})$. We let $D \subset F$ be a cell, and prove that $D \cap \{ a \} \times \mathbb{R}^{m}$ has cardinality at most $1$ for each $b$. From the definition, the cell $D$ is obtained as a sequence $D_{1}, D_{2}, \hdots, D_{n}, D_{n+1}, \hdots, D_{n+m}$ where $D_{i}$ is a cell in $\mathbb{R}^{i}$ and $D_{i+1}$ is obtained from $D_{i}$ in one of the two ways described in (ii) above. 

For each $1 \leq k \leq m$ we have a natural projection
\begin{equation}
\label{cellproj}
\pi_{k} : D_{n+k} \cap \{ a \} \times \mathbb{R}^{k} \to D_{n+(k-1)} \cap \{ a \} \times \mathbb{R}^{k-1} .
\end{equation}
For $k = 1$ this is just a map $D_{n+1} \cap \{ a \} \times \mathbb{R} \to \{ a \}$. We observe that all such maps are necessarily surjective, which is immediate from the two possibilities for the construction of $D_{n+k}$ from $D_{n+(k-1)}$ in (ii). Thus the hypothesis that $D_{n+m} \cap \{ a \} \times \mathbb{R}^{m}$ is finite in fact implies that $D_{n+k} \cap \{ a \} \times \mathbb{R}^{k} \to \{ a \}$ is finite for all $0 \leq k \leq m$, since each map $D_{n+m} \cap \{ a \} \times \mathbb{R}^{m} \to D_{n+k} \cap \{ a \} \times \mathbb{R}^{k}$ is surjective as it is a composition of surjective maps. 

It then suffices to prove inductively on $0 \leq k \leq m$ that $D_{n+k} \cap \{ a \} \times \mathbb{R}^{k}$ is a singleton, which at each stage amounts to proving that the (necessarily unique, by induction) fibre of $\pi_{k}$ is a singleton. Notice the finiteness of $D_{n+m} \cap \{ a \} \times \mathbb{R}^{k}$ means that $D_{n+k}$ cannot be obtained from $D_{n+(k-1)}$ by a construction of the second type in (ii) above, since then $D_{n+k} \cap \{ a \} \times \mathbb{R}^{k}$ would contain the infinite set
\[ \{ (x_0, y) : f(x_0) < y < g(x_0) \} \]
for the (necessarily unique, by induction) point $x_0 \in D_{n+(k-1)} \cap \{ a \} \times \mathbb{R}^{k-1}$. 

Hence $D_{n+k}$ is obtained from $D_{n+(k-1)}$ by a construction of the form
\[ D_{n+k} = \{ (x, y) \in D_{n+(k-1)} \times \mathbb{R} : y = f(x) \} . \]
But then the induction hypothesis that $D_{n+(k-1)} \cap \{ a \} \times \mathbb{R}^{k-1} = \{ x_0 \}$ is a singleton implies that $D_{n+k} \cap \{ a \} \times \mathbb{R}^{k} = \{ (x_{0}, f(x_0)) \}$, so the result follows.
\end{proof}

\begin{cor}
\label{finitepartcor}
Let $f : B \to A$ be a definable function with finite fibres. Then there exists a definable partition $B = B_{1} \sqcup \cdots \sqcup B_{k}$ such that $B_{i} \to A$ is injective.
\end{cor}

\begin{proof}
Let $F \subset B \times A$ be the graph of $f$. Then after swapping the order of the coordinates, we obtain from \autoref{finitepartlem} a definable partition $F = F_{1} \sqcup \cdots \sqcup F_{k}$ such that $F_{i} \cap B \times \{ a \}$ has cardinality at most $1$ for each $i$ and $a$. Define $B_{i}$ to be the projection of $F_{i}$ and observe that $F_{i}$ is just the graph of $B_{i} \to A$. 
\end{proof}

\subsection{Stripes}

\subsubsection{Parametrization of Stripes}
\label{paramstripesec}

We work in the setup of \S\ref{siegesetsec}. We additionally fix a fibre $V = \mathbb{V}_{s_{0}}$, the choice of which is unimportant, and identify $(G,D) = (\mathbf{G}_{S}, D_{S})$ with their realizations at $s_{0}$; in particular, $G$ is a subgroup of $\GL(V)$ which lies in the group $\textrm{SAut}(V, Q)$ of linear mappings which preserve the polarizing form $Q$ on $V$ and have determinant one. Fix a Hodge-theoretic Siegel set $\mathfrak{O} = \Omega A_{t} K_{o} \cdot o \subset D$ associated to $(P, o, t)$, and set $\mathfrak{G} = \Omega A_{t} K_{o}$. We will denote the natural orbit map $G(\mathbb{R}) \to D$ given by $g \mapsto g \cdot o$ by $q$. In what follows all definable sets are regarded as definable relative to the structure $\mathbb{R}_{\textrm{an,exp}}$. We let $h \in D$ be a point with Mumford-Tate group $M_{h} \subset G$, and fix another point $o \in D_{M} = M_{h}(\mathbb{R}) \cdot h$. Then $o$ has Mumford-Tate group contained in $M := M_{h}$. We write $\ch{D}$ for the compact dual of $D$, which is a projective algebraic variety containing $D$ as an open submanifold. It can be identified with the orbit of $G(\mathbb{C})$ on $F^{\bullet}_{o}$ in the space of all Hodge flags on $V$, where $F^{\bullet}_{o}$ is the Hodge flag corresponding to $o$. 

\begin{notn}
Suppose $X$ is a $\mathbb{Q}$-algebraic variety and $\mathfrak{S} \subset X(\mathbb{R})$ is a subset. We write $\mathfrak{S}(\mathbb{Q})$ for the intersection $\mathfrak{S} \cap X(\mathbb{Q})$.
\end{notn}

\begin{notn}
Given a group $G$ with subgroup $H$, we write $\mathcal{N}_{G}(H)$ (resp. $\mathcal{C}_{G}(H)$) for the normalizer (resp. centralizer) in $G$ of $H$. If both $G$ and $H$ are algebraic groups, we interpret $\mathcal{N}_{G}(H)$ (resp. $\mathcal{C}_{G}(H)$) as an algebraic group.
\end{notn}

\begin{notn}
For an algebraic group $G$, we write $G^{\circ}$ for its identity component.
\end{notn}

\begin{rem}
\label{centrem}
If $H$ is a connected reductive subgroup of a connected reductive group $G$, then $H \cdot \mathcal{C}_{G}(H)^{\circ} = \mathcal{N}_{G}(H)^{\circ}$ (c.f. \cite{27349}). This is often useful for understanding $\mathcal{N}_{\mathbf{G}_{S}}(\mathbf{M})$ when $\mathbf{M}$ is a Mumford-Tate subgroup of $\mathbf{G}_{S}$. 
\end{rem}

We set $D_{M} = M(\mathbb{R}) \cdot o$ to be the $M(\mathbb{R})$-orbit of $o$ in $D$, and $\ch{D}_{M} = M(\mathbb{C}) \cdot F^{\bullet}_{o} \subset \ch{D}$ to be the corresponding orbit in $\ch{D}$, where $F^{\bullet}_{o}$ is the Hodge flag corresponding to $o$. We write $\ch{\NL}_{M} \subset \ch{D}$ for the $\mathbb{Q}$-algebraic locus considering all Hodge flags whose set of Hodge tensors contains those tensors fixed by $M$ (c.f. \cite[Ch. 2]{GGK}). 

\begin{prop}
\label{FDMprop}
Let $\mathcal{F}(D_{M})$ be the Zariski closure in $G$ of
\[ \{ g \in G(\mathbb{R}) : g D_{M} \subset D_{M} \} , \]
and $\mathcal{F}(\ch{D}_{M})$ the algebraic group defined by 
\[ \mathcal{F}(\ch{D}_{M})(\mathbb{C}) := \{ g \in G(\mathbb{C}) : g \ch{D}_{M} \subset \ch{D}_{M} \} . \]
Then if $I$ is the identity component of $\mathcal{N}_{G}(M)$, we have
\begin{equation}
\label{canFDMinc}
I \subset \mathcal{F}(D_{M}) \subset \mathcal{F}(\ch{D}_{M}) .
\end{equation}
Moreover, the identity components of all three groups agree.
\end{prop}

\begin{proof}
We first observe that $g D_{M} \subset D_{M}$ implies $g \ch{D}_{M} \subset \ch{D}_{M}$ because $D_{M}$ is open in the irreducible variety $\ch{D}_{M}$. This implies the latter inclusion $\mathcal{F}(D_{M}) \subset \mathcal{F}(\ch{D}_{M})$. For the former we use \cite[VI.A.3]{GGK} to obtain that $\mathcal{N}_{G}(M) \ch{D}_{M} \subset \ch{\NL}_{M}$, and note that by \cite[VI.B.1]{GGK} $\ch{D}_{M}$ is a (geometrically) connected component of $\ch{\NL}_{M}$. As a consequence of \cite[VI.B.1]{GGK}, the (geometrically) irreducible components of $\ch{\NL}_{M}$ agree with the (geometrically) connected components. This implies that geometrically connected group $I$ preserves $\ch{D}_{M}$. Then if $g \in I(\mathbb{R})$ one necessarily has that $g D_{M} \subset D \cap \ch{D}_{M} = D_{M}$ (using \cite[VI.B.11]{GGK}), hence $I(\mathbb{R}) \subset \mathcal{F}(D_{M})(\mathbb{R})$. But $I(\mathbb{R})$ is Zariski dense in the connected $\mathbb{Q}$-group $I$ since connected algebraic groups over a field of characteristic zero are unirational \cite[Theorem 18.2]{lagbor}. 

Now for the claim about identity components. We start by using unirationality to observe that $\mathcal{F}(D_{M})(\mathbb{Q})^{\circ} := \mathcal{F}(D_{M})(\mathbb{Q}) \cap \mathcal{F}(D_{M})^{\circ}(\mathbb{R})$ is Zariski dense in $\mathcal{F}(D_{M})^{\circ}$. Now consider an element $g \in \mathcal{F}(D_{M})(\mathbb{Q})^{\circ}$. We have $g D_{M} \subset D \cap \ch{D}_{M} = D_{M}$, and so $g$ sends a generic point of $D_{M}$ with Mumford-Tate group $M$ to another such point. Applying \cite[VI.A.3]{GGK} this tells us that $g \in \NL_{M}(\mathbb{R})$. Then necessarily $g \in I(\mathbb{R})$ since $\mathcal{F}(D_{M})^{\circ}$ is connected. The analogous argument works for $\mathcal{F}(\ch{D}_{M})$. 
\end{proof}

\begin{defn}
A $\mathbb{Q}$-subgroup $H \subset \mathcal{N}_{G}(M)$ is said to generate $D_{M}$ if $H(\mathbb{R}) \cdot o = D_{M}$.
\end{defn}

We set $V^{\otimes} = \bigoplus_{a,b \geq 0} V^{\otimes a} \otimes (V^{*})^{\otimes b}$. Each flag $F^{\bullet} \in \ch{D}$ induces a filtration on $V^{\otimes}$ which we denote with the same notation.

\begin{lem}
\label{Htransconstrlem}
Suppose that $H \subset \mathcal{N}_{G}(M)$ generates $D_{M}$ and consider a point $g H \in (G/H)(\mathbb{Q})$ such that $D_{g} = g \ch{D}_{M} \cap D$ is non-empty. Then there exists a Hodge structure $h_{g} \in D_{g}$ with Mumford-Tate group $M_{g} \subset g M g^{-1}$ such that $D_{g} = D_{M_{g}} := M_{g}(\mathbb{R}) \cdot h_{g}$. If $g H$ admits a representative $g \in G(\mathbb{R})$, one has $D_{g} = g D_{M}$ and $M_{g} = g M g^{-1}$. 
\end{lem}

\begin{proof}
Note that because $H \subset \mathcal{N}_{G}(M)$, the group $g M g^{-1}$ is independent of the representative $g$ chosen. Thus $g M g^{-1}$ is $\mathbb{Q}$-algebraic because $g H$ is. Let $g F^{\bullet}$ be a Hodge flag corresponding to a point $h_{g} \in D_{g}$, where $F^{\bullet}$ is a Hodge flag corresponding to a point of $\ch{D}_{M}$. Write $\textrm{Fix}(M) \subset V^{\otimes}_{\mathbb{C}}$ for the vector subspace of tensors fixed by $M$. Choosing $F^{0}$ sufficiently generic, this is exactly the $\mathbb{C}$-span of $F^{0} \cap V^{\otimes}_{\mathbb{Q}}$. One has $g  \textrm{Fix}(M) = \textrm{Fix}(g M g^{-1})$. As the group $g M g^{-1}$ is defined over $\mathbb{Q}$, the complex vector space $g \textrm{Fix}(M)$ has a natural underlying $\mathbb{Q}$-structure $[g \textrm{Fix}(M)]_{\mathbb{Q}}$, and by construction $[g \textrm{Fix}(M)]_{\mathbb{Q}} \subset g F^{0} \cap V^{\otimes}_{\mathbb{Q}}$. Then the Mumford-Tate group $M_{g}$ of $h_{g}$ fixes all tensors in $g F^{0} \cap V^{\otimes}_{\mathbb{Q}}$, hence lies inside $g M g^{-1}$. Taking $h_{g}$ sufficiently generic and using that $g \ch{D}_{M}$ is irreducible, the points of $D_{g}$ have Mumford-Tate group contained in a common $\mathbb{Q}$-algebraic Mumford-Tate group $M_{g} \subset g M g^{-1}$. Let $D_{M_{g}} = M_{g}(\mathbb{R}) \cdot h_{g}$. We then have
\[ D_{g} \subset D \cap \ch{\NL}_{M_{g}} = D_{M_{g}} , \]
where we have applied \cite[VI.B.11]{GGK}. 

Now $g \ch{\NL}_{M}$ is exactly the locus of Hodge flags whose set of Hodge tensors contains $[g \textrm{Fix}(M)]_{\mathbb{Q}}$. It thus follows that $\ch{\NL}_{M_{g}} \subset g \ch{\NL}_{M}$. For both the varieties $\ch{\NL}_{M_{g}}$ and $g \ch{\NL}_{M}$ their geometrically irreducible components agree with their geometrically connected components (this is a consequence of \cite[VI.B.1]{GGK}), so the inclusion $\ch{\NL}_{M_{g}} \subset g \ch{\NL}_{M}$ also implies an inclusion $\ch{D}_{M_{g}} \subset g \ch{D}_{M}$ of components passing through $h_{g}$. Then
\[ D_{M_{g}} = D \cap \ch{D}_{M_{g}} \subset D \cap g \ch{D}_{M} = D_{g} . \]

It remains to show that $M_{g} = g M g^{-1}$ and that $D_{g} = g D_{M}$ under the assumption $g \in G(\mathbb{R})$. For the second claim, observe that 
\[ g^{-1} D_{g} = g^{-1} (g \ch{D}_{M} \cap D) = \ch{D}_{M} \cap D = D_{M} , \]
where we again apply \cite[VI.B.11]{GGK}. Hence $g D_{M} = D_{g}$. Finally $g^{-1} M_{g} g$ is the generic Mumford-Tate group of a point in $g^{-1} D_{g} = D_{M}$, hence $g^{-1} M_{g} g = M$.
\end{proof}

\begin{defn}
A Mumford-Tate domain $D_{g}$ obtained as in \autoref{Htransconstrlem} is called an $H$-translate of $D_{M}$.
\end{defn}

We thus obtain a map 
\[ \eta : (G/H)(\mathbb{Q}) \to \{ H\textrm{-translates of }D_{M} \} \]
which sends $g H$ to $D_{M_{g}}$. We write $(G/H)(\mathbb{R})^{\circ} \subset (G/H)(\mathbb{R})$ for the connected component containing $H$. We then write $\eta^{\circ}$ for the restriction of $\eta$ to $(G/H)(\mathbb{Q})^{\circ} = (G/H)(\mathbb{Q}) \cap (G/H)(\mathbb{R})^{\circ}$. 

\begin{lem}
\label{fibreetaliesinfib}
Each fibre of $\eta$ lies inside a fibre of $(G/H)(\mathbb{C}) \to (G/\mathcal{F}(\ch{D}_{M}))(\mathbb{C})$. Each fibre of $\eta^{\circ}$ lies inside a fibre of $G(\mathbb{R})/H(\mathbb{R}) \to G(\mathbb{R})/\mathcal{F}(D_{M})(\mathbb{R})$.
\end{lem}

\begin{proof}
If $\eta(g H) = \eta(g' H)$ then $g \ch{D}_{M} = g' \ch{D}_{M}$ which implies that $g'^{-1} g \in (\mathcal{F}(\ch{D}_{M}))(\mathbb{C})$. Similarly if $\eta^{\circ}(g H) = \eta^{\circ}(g' H)$ then $g D_{M} = g' D_{M}$ and hence $g'^{-1} g \in \mathcal{F}(D_{M})(\mathbb{R})$. 
\end{proof}

\begin{defn}
Suppose that $\mathfrak{T} \subset \mathfrak{O}$ is a subset containing $o$, and that $H \subset G$ is a $\mathbb{Q}$-algebraic subgroup which generates $D_{M}$. By a stripe of $H$ in $\mathfrak{T}$ we mean a non-empty intersection $\eta(g H) \cap \mathfrak{T}$ for some $g H \in (G/H)(\mathbb{Q})$. A stripe is said to be principal if we can choose $g H \in (G/H)(\mathbb{Q}) \cap (G/H)(\mathbb{R})^{\circ}$. 
\end{defn}

\begin{defn}
\label{genericstripedef}
Given a subset $\mathfrak{T} \subset \mathfrak{O}$ containing $o$ and a $\mathbb{Q}$-subgroup $H \subset G$ generating $D_{M}$, we say that a stripe 
\[ \mathfrak{H}_{g} := \eta(g H) \cap \mathfrak{T} = D_{M_{g}} \cap \mathfrak{T} \]
of $H$ in $\mathfrak{T}$ is generic if there is a point of $\mathfrak{H}_{g}$ with Mumford-Tate group $M_{g}$. We say it is uniformly generic if every connected component of $\mathfrak{H}_{g}$ has such a point.
\end{defn}

\begin{lem}
\label{genericwelldeflem}
In the context of \autoref{genericstripedef}, the genericity of the stripe $\mathfrak{H}_{g}$ does not depend on $g$; if $\mathfrak{H}_{g} = \mathfrak{H}_{g'}$ and there exists a point $t \in \mathfrak{H}_{g}$ with Mumford-Tate group $M_{g}$, then in fact $M_{g} = M_{g'}$ and $D_{M_{g}} = D_{M_{g'}}$. 
\end{lem}

\begin{proof}
Let $t \in \mathfrak{H}_{g} \subset D_{M_{g}}$ be a point with Mumford-Tate group equal to $M_{g}$. Then $t \in \mathfrak{H}_{g'} \subset D_{M_{g'}}$ so $M_{g} \subset M_{g'}$. On the other hand $\dim D_{M_{g}} = \dim D_{M} = \dim D_{M_{g'}}$, so $M_{g}(\mathbb{R}) \cdot t$ is open in $D_{M_{g'}}$. Since any open subset of $D_{M_{g'}}$ contains a Hodge-generic point it follows that $M_{g} = M_{g'}$. Likewise, $D_{M_{g}} = D_{M_{g'}}$. 
\end{proof}

We now obtain a diagram 
\begin{equation}
\label{stripemapdiag}
\begin{tikzcd}
 & \{ \textrm{generic stripes of }H\textrm{ in }\mathfrak{T} \} \arrow[d,hook,"\beta = \beta_{\mathfrak{T}}"] \\
(G/H)(\mathbb{Q}) \arrow[r,"\eta"] & \{ H\textrm{-translates of }D_{M} \} ,
\end{tikzcd}
\end{equation}
where the map $\beta$ sends a stripe $\mathfrak{H}_{g}$ to the (necessarily unique by \autoref{genericwelldeflem}) $H$-translate of $D_{M}$ which induces it. We wish to understand the inverse image $\eta^{-1}(\textrm{im}(\beta))$ by constructing, using $\mathfrak{T}$, a subset of $(G/H)(\mathbb{R})$ which contains it. For this we define
\begin{align*}
\mathfrak{Z} := \textrm{im}[q^{-1}(\mathfrak{T}) \to (G/H)(\mathbb{R})] \cap (G/H)(\mathbb{R})^{\circ} \subset (G/H)(\mathbb{R}) .
\end{align*}

\begin{lem}
The restriction of $\eta$ to $\mathfrak{Z}(\mathbb{Q})$ induces a surjection onto the principal stripes of $H$ in $\mathfrak{T}$.
\end{lem}

\begin{proof}
Suppose that $g H \in (G/H)(\mathbb{Q}) \cap (G/H)(\mathbb{R})^{\circ}$ is such that $\eta(gH) \cap \mathfrak{T}$ is non-empty. Since $g H \in (G/H)(\mathbb{R})^{\circ} \subset G(\mathbb{R}) \cdot H$, we may choose $g \in G(\mathbb{R})$. From \autoref{Htransconstrlem} we then know that $g D_{M} = D_{g} = \eta(g H)$, so $g D_{M} \cap \mathfrak{T}$ is non-empty. Since $D_{M} = H(\mathbb{R}) \cdot o$, we may therefore choose $h \in H(\mathbb{R})$ such that $g h \cdot o = t \in \mathfrak{T}$. Then $g h \in q^{-1}(\mathfrak{T})$ and $g H \in \mathfrak{Z}$.
\end{proof}

\subsection{Point Counting in Coset Spaces}
\label{pointcountingcosetsec}

In this section we explain how to count rational points inside $(G/H)(\mathbb{R})$, at least for appropriately chosen $G$ and $H$. For this we will use some results of \cite{zbMATH05902967}. 

\subsubsection{Measures}

\begin{defn}
We say an algebraic group is scss if it is simply-connected and semisimple.
\end{defn}

Suppose that $H \subset G$ is an inclusion of $\mathbb{Q}$-algebraic scss groups and write $G/H$ for the quotient $\mathbb{Q}$-variety. Given Haar measures $\mu_{H}$ and $\mu_{G}$ on $H(\mathbb{A})$ and $G(\mathbb{A})$, respectively, we say that a measure $\mu_{G/H} := \mu_{Z}$ on $Z := G(\mathbb{A})/H(\mathbb{A})$ is compatible with $\mu_{H}$ and $\mu_{G}$ if for any compactly supported function $f$ on $G(\mathbb{A})$ we have

\[ \int_{G(\mathbb{A})} f(t)\,  \mu_{G}(t) = \int_{Z} \int_{H(\mathbb{A})} f(gu) \, d\mu_{H}(u) \, d\mu_{Z}(gH(\mathbb{A})) . \]

\noindent In our case, because both $G$ and $H$ are unimodular, such a measure exists by \cite[Thm. 5.2.1, 5.2.2]{zbMATH06960207} and is uniquely determined up to a scalar. This scalar can be fixed by requiring that $\mu_{H}$ (resp. $\mu_{G}$) induces a probability measure on the quotient $H(\mathbb{A})/H(\mathbb{Q})$ (resp. $G(\mathbb{A})/G(\mathbb{Q})$); note that these quotients have finite volume as a consequence of \cite[Thm. 1]{zbMATH03194969} \cite{zbMATH03194968}. We will always adopt the convention whenever considering a triple $(G, H, Z = G(\mathbb{A})/H(\mathbb{A}))$ with $G$ and $H$ scss that the measures $\mu_{G}, \mu_{H}$ and $\mu_{G/H} = \mu_{Z}$ have been chosen in this way. 

In a situation where we have a representation $\rho : G \to \GL(V)$ on a $\mathbb{Q}$-vector space $V$ such that $H$ is the stabilizer in $G$ of some $v \in V$, we will also denote by $\mu_{G/H}$ the induced measure on the adelic points of the orbit variety $G \cdot v \simeq G/H$. Using the inclusion $\mathbb{Q}_{p} \hookrightarrow \mathbb{A}$ we also obtain induced measures $\mu_{G/H,p}$ on $(G/H)(\mathbb{Q}_{p})$ for each prime $p$, and likewise a measure $\mu_{G/H,\infty}$ at the infinite place.

\subsubsection{Orbit Asymptotics}
\label{orbitasympsec}

We now recall a theorem proven in \cite{zbMATH05902967}. We note that \cite{zbMATH05902967} uses right-coset spaces instead of left ones, the latter being our convention. The notation $C_{c}(-)$ refers to compactly supported continuous functions on the topological space inside the brackets. In what follows we set $\mu = \mu_{Z}$.

\begin{prop}[Prop. 5.3 in \cite{zbMATH05902967}]
\label{EOkeyprop}
Suppose that both $G$ and $H$ are scss $\mathbb{Q}$-groups with $H \subset G$ a maximal connected $\mathbb{Q}$-subgroup. Set $Z = G(\mathbb{A}) / H(\mathbb{A})$. 
Then for any well-rounded sequence $\{ B_{T} \subset Z \}$ of compact subsets whose volume diverges as $T \to \infty$, we have,
\[ \# (G/H)(\mathbb{Q})^{\circ} \cap B_{T} = \# (G(\mathbb{Q}) [1] \cap B_{T}) \sim \mu(B_{T}) , \]
where $[1] \in Z$ is the class of the identity.
\end{prop}

\vspace{0.5em}

The statement of the proposition uses the following definition:

\begin{defn}
A family $B_{T} \subset Z$ of compact subsets is called well-rounded if there exists $c > 0$ such that for every small $\ep > 0$, there exists a neighbourhood $U_{\ep}$ of $1$ in $G(\mathbb{A})$ such that for all sufficiently large $T$,
\[ (1 - c \cdot \ep) \mu(U_{\ep} B_{T}) \leq \mu(B_{T}) \leq (1 + c \cdot \ep) \mu( \cap_{u \in U_{\ep}} u B_{T} ) . \]
\end{defn}

\vspace{0.5em}

We note that our definition of ``well-rounded'' is the special case of the notion ``$W$-well-rounded'' appearing in \cite[Def. 5.1]{zbMATH05902967} with $W = 1$.

\begin{proof}[Proof of (\ref{EOkeyprop}):]
We start by focusing on the asymptotic claim. This is a special case of \cite[Prop. 5.3]{zbMATH05902967}. Note that the authors of \cite{zbMATH05902967} use $G$ where we use $G(\mathbb{A})$, $L$ where we use $H(\mathbb{A})$, and we take $\Gamma = G(\mathbb{Q})$. Here we have used the normalization $\mu_{G}(G(\mathbb{A}) / G(\mathbb{Q})) = 1$ as well as $\mu_{H}(H(\mathbb{A}) / H(\mathbb{Q})) = 1$. The equidistribution hypothesis of \cite{zbMATH05902967} (i.e., the convergence of the expression involving integrals) is a consequence of \cite[Cor. 4.14]{zbMATH05902967}, where we note that $G_{W} = G(\mathbb{A})$ for us.

For the first equality we note that, because $H$ is simply connected, there is exactly one $G(\mathbb{Q})$ orbit in each $G(\mathbb{R})$ orbit inside $G/H$ (see the proof of \cite[Cor. 1.9]{zbMATH05902967}). Thus it makes sense to identify $(G/H)(\mathbb{Q})^{\circ}$ with $G(\mathbb{Q}) \cdot H(\mathbb{R}) \subset G(\mathbb{R})/H(\mathbb{R})$, and the result follows. 
\end{proof}

To apply \autoref{EOkeyprop}, it suffices to construct a sequence of well-rounded sets. To do this we proceed as follows, following \cite[Pf. of Cor. 1.9]{zbMATH05902967}. We let $\rho : G \to \GL_{m}$ be a representation of $G$ satisfying the condition that $H \subset G$ is identified with the stabilizer of some vector $v \in \mathbb{Q}^{m}$. Fix a compact measurable subset $\Omega \subset G(\mathbb{R}) \cdot v$ with boundary measure zero and positive volume. We then consider, for each positive integer $\ell$, the sets
\begin{equation}
\label{wellroundsets}
B_{\ell} := \{ (x_{p}) \in G(\mathbb{A}) \cdot v : x_{\infty} \in \Omega, \ \| x_{p} \|_{p} \leq p^{e_{p}} \textrm{ for }p\textrm{ a finite prime} \}, \hspace{1.5em} \ell = \prod_{p} p^{e_{p}} .
\end{equation}
We observe that $B_{\ell}$ is well-rounded with $c = 1$. Indeed, we can consider the subgroup $\prod_{p} G(\mathbb{Z}_{p}) \subset G(\mathbb{A})$, which we may observe preserves $B_{\ell}$. Then taking a neighbourhood of $1$ of the form $U_{\ep} = K_{\ep} \times \prod_{p} G(\mathbb{Z}_{p})$ for $K_{\ep}$ a sufficiently small compact neighbourhood of $1 \in \mathbb{R}$ depending on $\ep$ we easily see that $B_{\ell}$ is well-rounded. From the calculation in \cite[Cor. 7.7]{zbMATH06130496} (c.f. the proof of Cor. 1.9 in \cite{zbMATH05902967}), one also sees that $\mu(B_{\ell}) \to \infty$. We conclude from \autoref{EOkeyprop} that
\begin{cor}
\label{pointcountcor}
In the above setup,
\begin{equation}
\label{pointcountconverrges}
\# (G/H)(\mathbb{Q})^{\circ} \cap B_{\ell} \sim \mu(B_{\ell}) .
\end{equation}
\end{cor}

We remark that in the argument appearing in \cite[Pf. of Cor. 1.9]{zbMATH05902967} the representation $\rho : G \to \GL_{m}$ is assumed to be faithful, but as we have just seen for the purpose of establishing \autoref{pointcountcor} this is not need. Supposing now that we define

\begin{equation}
\label{wellroundsets}
B^{\circ}_{\ell} := \{ (x_{p}) \in G(\mathbb{A}) \cdot v : x_{\infty} \in \Omega^{\circ}, \ \| x_{p} \|_{p} \leq p^{e_{p}} \textrm{ for }p\textrm{ a finite prime} \}, \hspace{1.5em} \ell = \prod_{p} p^{e_{p}} 
\end{equation}
where $\Omega^{\circ} \subset \Omega$ denotes the interior of $\Omega$, we also have

\begin{cor}
\label{pointcountingininterior}
In the above setup, 
\begin{equation}
\label{pointcountconverrges}
\# (G/H)(\mathbb{Q})^{\circ} \cap B^{\circ}_{\ell} \sim \mu(B^{\circ}_{\ell}) = \mu(B_{\ell}) .
\end{equation}
\end{cor}

\begin{proof}
Using the structure of the product measure, we have $\mu(B_{\ell}) = \mu(\Omega) \cdot \mu_{f}(B'_{\ell})$, where $B'_{\ell} \subset G(\mathbb{A}_{f}) \cdot v$ is defined as $B_{\ell}$ is except with the condition $x_{\infty} \in \Omega^{\circ}$ at the real place. We may approximate $\Omega^{\circ}$ from below as a union of compact subsets $\Omega_{j}$, and by applying \autoref{EOkeyprop} to a sequence of sets $B_{\ell j}$ constructed with the $\Omega_{j}$ we learn that $\# (G/H)(\mathbb{Q})^{\circ}  \cap \Omega^{\circ}$ is asymptotically greater than $\mu(\Omega_{j}) \cdot \mu_{f}(B'_{\ell})$ for each $j$. Since $\mu(\Omega_{j}) \to \mu(\Omega)$ we obtain the result. 
\end{proof}

\section{Proofs of the General Theorems}

We define a subset $\mathcal{S}(\mathbb{C})$ of points in the tensorial Hodge locus which will be of interest.

\[ \mathcal{S} = \left\{ s \in S(\mathbb{C}) : \begin{matrix}
(\mathbf{G}_{s}, D_{s}) \subset (\mathbf{M}', D')\textrm{ where }(\mathbf{M}', D')\textrm{ is a } \\ \textrm{ Hodge subdatum }\mathbb{Q}\textrm{-coset equivalent to }(\mathbf{M}, D_{M}) \end{matrix} \right\} . \]

\begin{defn}
\label{suffgenhypplanedef}
We say a hyperplane section $L$ of $S$ of codimension $d$ is sufficiently general if all of the following conditions hold:
\begin{itemize}
\item[(i)] $L$ is smooth and irreducible;
\item[(ii)] $L$ intersects every Hodge locus component of dimension $d$ in a reduced set of points of cardinality equal to its degree;
\item[(iii)] $L$ does not intersect any (possibly tensorial) Hodge locus component of dimension smaller than $d$;
\item[(iv)] $L$ is not contained in any component of the tensorial Hodge locus, and the algebraic monodromy group of $\mathbb{V}$ agrees with that of $\restr{\mathbb{V}}{S \cap L}$; and
\item[(v)] the restriction of the period map $\varphi$ to $S \cap L$ is quasi-finite.
\end{itemize}
\end{defn}

\noindent Note that sufficiently general hyperplanes sections exist assuming $d \geq \dim S - \dim \varphi(S)$. To explain why one can achieve (v), we note that the map $\varphi$ factors as $\varphi : S \to T \hookrightarrow \Gamma \backslash D_{S}$ with $S \to T$ algebraic (see \cite{zbMATH07662555}), and all the components of the (tensorial) Hodge locus are pulled back from algebraic subvarieties of $T$. In particular to achieve (v) it is enough to require that $L$ does not intersect any fibre of $S \to T$ in a positive dimensional locus, which is true away from a closed locus in the parameter space of hyperplane sections assuming $d \geq \dim S - \dim \varphi(S)$. 

\subsection{Proof of \autoref{mainthmupbound}}

We will fix a $d$, and prove the same statement but with $Z(v')$ replaced by $Z(v')_{d}$. This suffices, since one can always obtain (\ref{degpointcorres}) by summing over such inequalities for all possible $d$ and combining the constants. 

\vspace{0.5em}

\noindent \textbf{Reduction to Local Point Counting:} Now letting $L \subset S$ be a sufficiently general hyperplane section of codimension $d$, we can replace $S$ with $S \cap L$, and reduce to proving the same result for the restricted variation on the new space. Note that condition (iii) in particular implies that the Mumford-Tate groups of the Hodge structures above the points of $L \cap Z(v')$ are equal to the Mumford-Tate group of the corresponding component of $Z(v')$. After all these changes one can take $d = 0$ and reduce to estimating the degrees of the union of components $Z(v')_{0}$.

\vspace{0.5em}

We now use \autoref{landsinsiegelset} to take an open cover $S = \bigcup_{i = 1}^{n} B_{i}$ such that each $B_{i}$ admits a definable local period map $\psi_{i} : B_{i} \to D_{S}$ landing inside a Siegel set $\mathfrak{O}_{i} \subset D_{S}$, with $\mathfrak{O}_{i} = \Omega_{i} A_{t,i}  K_{o_{i}} \cdot o_{i}$, with $\mathcal{S}$ the set defined above matching the one in \autoref{landsinsiegelset}. We set $\mathfrak{G}_{i} = \Omega_{i} A_{t,i}  K_{o_{i}}$. We therefore have, for each $B_{i}$ which intersects $\mathcal{S}$, a point $s_{i} \in B_{i}$ such that the Hodge datum of $o_{i} := \psi_{i}(s_{i})$ is contained in a Hodge datum $\mathbb{Q}$-coset equivalent to $(\mathbf{M}, D_{M})$. We may moreover use \autoref{finitepartcor} and the fact that $\varphi$ is quasi-finite to ensure that each map $\psi_{i}$ is actually injective (we no longer require that $o_{i}$ is the image of some point $s_{i} \in B_{i}$, just that the Hodge datum of $o_{i}$ be contained in a Hodge datum $\mathbb{Q}$-coset equivalent to $(\mathbf{M}, D_{M})$).

By shrinking the $B_{i}$ to closed subsets $C_{i} \subset B_{i}$, we may assume that $S = \bigcup_{i=1}^{n} C_{i}$ but that the sets $\{ C_{1}, \hdots, C_{n} \}$ have disjoint interiors, which we denote by $C^{\circ}_{i}$. Moreover, we can choose the $C_{i}$ so that the boundary $C_{i} \setminus C^{\circ}_{i}$ avoids the countably many points $\bigcup_{v'} Z(v')_{0}$ for each $i$. We then have that
\begin{equation}
\label{reltopointcountingeq}
\left[ \sum_{\substack{v' \in O(\mathbb{Q})^{\circ} \\ \nu(v') | \nu}} \deg Z(v')_{0} \right] = \sum_{i = 1}^{n} \left| C^{\circ}_{i} \cap \left(\bigcup_{\substack{v' \in O(\mathbb{Q})^{\circ} \\ \nu(v') | \nu}} Z(v')_{0} \right) \right| .
\end{equation}
It will therefore suffice to estimate the asymptotic size as $\nu \to \infty$ of each of the summands on the right. 

\vspace{0.5em}

\noindent \textbf{Relating to Counting Stripes: } We now let $q_{i} : \mathbf{G}_{S}(\mathbb{R}) \to D_{S}$ be the natural orbit map $g \mapsto g \cdot o_{i}$. By our choice above, the Hodge structure $o_{i}$ has Hodge datum contained in a datum $\mathbb{Q}$-coset equivalent to $(\mathbf{M}, D_{M})$ we denote by $(\mathbf{M}_{i}, D_{M_{i}})$. Because $(\mathbf{M}_{i}, D_{M_{i}})$ is $\mathbb{Q}$-coset equivalent to $(\mathbf{M}, D_{M})$, there is $g \in \mathbf{G}_{S}(\mathbb{R})$ relating the two Hodge subdatums such that $g \mathbf{M}$ is defined over $\mathbb{Q}$. Using the isomorphism $\mathbf{G}_{S}/\mathbf{M} \cong O$, this means that $v_{i} := (g \mathbf{M}) \cdot v$ is a $\mathbb{Q}$-vector and $\mathbf{M}_{i}$ is exactly the stabilizer in $\mathbf{G}_{S}$ of $v_{i}$. Note this means that $v_{i}$ is Hodge for $o_{i}$. 

Let $r_{i} : \mathbf{G}_{S}(\mathbb{R}) \to (\mathbf{G}_{S}/\mathbf{M}_{i})(\mathbb{R})$ be the projection. We set $\mathfrak{T}_{i} = \psi_{i}(C^{\circ}_{i})$ and $\mathfrak{Z}_{i} = r_{i}(q^{-1}_{i}(\mathfrak{T}_{i}))$. Using the injectivity of $\psi_{i}$ we may then compute the summands in (\ref{reltopointcountingeq}) by
\begin{equation}
\label{pointcountinperdom}
\left| C^{\circ}_{i} \cap \left(\bigcup_{\substack{v' \in O(\mathbb{Q})^{\circ} \\ \nu(v') | \nu}} Z(v')_{0} \right) \right| = \left| \bigcup_{\substack{g \in (\mathbf{G}_{S}/\mathbf{M}_{i})(\mathbb{Q})^{\circ} \\ \nu(g \cdot v_{i}) | \nu}} [\mathfrak{T}_{i} \cap g D_{M_{i}}]_{0} \right| .
\end{equation}

By construction, one also has
\begin{equation}
\mathfrak{Z}_{i} = \{ g \in (\mathbf{G}_{S}/\mathbf{M}_{i})(\mathbb{R})^{\circ} : (\mathfrak{T}_{i} \cap g D_{M_{i}}) \neq \varnothing \} .
\end{equation}
Because $\mathfrak{T}_{i}$ is definable and the sets $g D_{M_{i}}$ lie in a common definable family, the number of isolated points in $\mathfrak{T}_{i} \cap g D_{M_{i}}$ is uniformly bounded by some constant $\kappa$ independent of $g$. This implies that
\begin{align*}
\left| \bigcup_{\substack{g \in (\mathbf{G}_{S}/\mathbf{M}_{i})(\mathbb{Q})^{\circ} \\ \nu(g \cdot v_{i}) | \nu}} [\mathfrak{T}_{i} \cap g D_{M_{i}}]_{0} \right| &\leq \kappa \cdot |\{ g \in \mathfrak{Z}_{i} \cap (\mathbf{G}_{S}/\mathbf{M}_{i})(\mathbb{Q})^{\circ} : \nu(g \cdot v) | \nu \}| \\
&\leq  \kappa \cdot |\{ v' \in (\mathfrak{G}_{i} \cdot v_{i})(\mathbb{Q}) : \nu(v') | \nu \}| ,
\end{align*}
where we use that $q_{i}^{-1}(\mathfrak{T}_{i}) \subset \mathfrak{G}_{i}$ as well as the isomorphism $\mathbf{G}_{S}/\mathbf{M}_{i} \cong \mathbf{G}_{S} \cdot v_{i}$. This is what we wanted to show.

\subsection{Proof of \autoref{mainthmlowbound}}

The statement is clearly local, and it suffices to take $B$ relatively compact admitting a local period map $\psi : B \to D_{S}$ which lifts $\varphi$. The strongly $\mathbb{V}$-likely hypothesis implies that $\mathcal{S}$ is dense in $S(\mathbb{C})$ by \cite[Thm. 3.5]{arXiv:2303.16179} (in fact even just $\mathbb{V}$-likely suffices here). Thus $B$ contains a point $s \in \mathcal{S}$. We set $o = \psi(s)$. After relabelling we may assume that the Hodge datum at $o$ is contained in $(\mathbf{M}, D_{M})$. We set $q : \mathbf{G}_{S}(\mathbb{R}) \to D_{S}$ to be the orbit map $g \mapsto g \cdot o$. 

We set $\mathcal{I} = \psi(B)$, and follow the proof given in \cite[\S4.4]{arXiv:2303.16179}. The proof given there explains that we have an open neighbourhood in $\mathbf{G}_{S}(\mathbb{R}) \times \mathcal{I}$ such that for each $(g, x)$ in this neighbourhood the intersection $g D_{M} \cap \mathcal{I}$ has dimension $d = \dim \mathcal{I} + \dim D_{M} - \dim D_{S}$ at $x$. Let us be more precise about what this neighbourhood looks like. One can start by defining $\mathcal{V} \subset \mathbf{G}_{S}(\mathbb{R})$ as the locus of $g \in \mathbf{G}_{S}(\mathbb{R})$ for which $g D_{M} \cap \mathcal{I}$ is non-empty. Then the argument in \cite[\S4.3]{arXiv:2303.16179} shows that the condition that $\dim_{x} (g D_{M} \cap \mathcal{I}) = d$ is open on $(g,x) \in \mathcal{V} \times \mathcal{I}$. The proof of \cite[\S4.4]{arXiv:2303.16179} then shows that, after removing a closed locus $\mathcal{C}$ of smaller dimension from $\mathcal{V} \times \mathcal{I}$, the intersection germs $(g D_{M} \cap \mathcal{I}, x)$ all have pure dimension $d$ and do not lie inside any translates of period subdomains of $D_{S}$ of smaller dimension. In particular, whenever $(g, x)$ is outside of $\mathcal{C}$ and $g D_{M}$ is a Mumford-Tate domain, the intersection $g D_{M} \cap \mathcal{I}$ is a generic stripe of $H = \mathbf{M}$ in $\mathcal{I}$, where we use the language of \S\ref{paramstripesec}. Using the definition of the product topology we can find a product $\mathcal{V}' \times \mathcal{I}' \subset (\mathcal{V} \times \mathcal{I}) \setminus \mathcal{C}$. Replacing $B$ with $B' = \psi^{-1}(\mathcal{I}')$ we may then assume that $\mathcal{C} = \varnothing$. 

Now choose $L$ sufficiently general. This means in particular that, for each of the intersections $L \cap Z(v')$, each point of $L \cap Z(v')$ has the same Mumford-Tate group as the component of $Z(v')$ in which it lies. In particular, each intersection $g D_{M} \cap \psi(B \cap L)$ with $(g \mathbf{M} g^{-1}, g D_{M})$ a Hodge datum $\mathbb{Q}$-coset equivalent to $(\mathbf{M}, D_{M})$ and $g \in \mathcal{V}$ is a generic stripe of $D_{M}$ in $\psi(B \cap L)$. Moreover the condition that $g D_{M}$ intersect $\psi(B \cap L)$ is still open on $\mathcal{V} \times \psi(B \cap L)$ (the analytic varieties $\psi(B \cap L)$ and $g D_{M}$ have complementary dimension inside $D_{S}$), and so after possibly shrinking $\mathcal{V}$ we can ensure that for each $g \in \mathcal{V}$ the intersection $g D_{M} \cap \psi(B \cap L)$ is a generic stripe of $D_{M}$ in $\mathfrak{T} := \psi(B \cap L)$.

Now let $U$ be the image of $\mathcal{V}$ in $O(\mathbb{R})^{\circ} \cong \mathbf{G}_{S}(\mathbb{R})/\mathbf{M}(\mathbb{R})$, where for the isomorphism we take the natural orbit map. Then the set $U$ is open inside the set $\mathfrak{Z}$ of \S\ref{paramstripesec}. We then consider the map 
\[ U(\mathbb{Q}) \xrightarrow{\eta^{\circ}} \{ \textrm{generic stripes of }H\textrm{ in }\mathfrak{T} \} \]
induced by the restriction of the map $\eta^{\circ}$ of \S\ref{paramstripesec}. Combining \autoref{fibreetaliesinfib}, \autoref{FDMprop} and the assumption that $\mathbf{M}$ has finite index in its normalizer, we learn that $\eta^{\circ}$ is quasi-finite with fibres of uniformly bounded size. Identifying $\mathbf{G}_{S}/\mathbf{M}$ with $O$ using the natural representation, we get that
\begin{equation}
\label{ptcountasymp}
\left| U \cap \frac{1}{\nu} \mathbb{V} \right| \asymp \sum_{\substack{g \in U(\mathbb{Q}) \\ \nu(g \cdot v) | \nu}} \left| [\psi(B \cap L) \cap g D_{M}]_{0} \right| . 
\end{equation}
Using the genericity of $L$, the sum on the right is then a lower bound for the sum in the statement. Since the symbol $\asymp$ implies the same expression with the symbol $\lesssim$ up to a positive real constant, the result follows.

\subsection{Proof of \autoref{Siegelorbitboundprop} and \autoref{NLlocusbound}}

\autoref{Siegelorbitboundprop} was shown in \S\ref{pointcountingsec}. Then \autoref{NLlocusbound} follows directly by combining \autoref{Siegelorbitboundprop} with \autoref{mainthmupbound}.

\subsection{Proof of \autoref{NLqbound}}

Given a Hodge structure $h$ on a lattice $V$, a vector $v \in V$ is Hodge for $h$ if and only if $v \otimes v$ is Hodge for the Hodge structure $h \otimes h$ on $V \otimes V$. (Indeed, the Mumford-Tate group $M_{h}$ of $h$ fixes $v$ if and only if it fixes $v \otimes v$ acting diagonally, and $M_{h \otimes h}$ is identified with the image of $M_{h}$ under the diagonal action since it is the $\mathbb{Q}$-Zariski closure of the diagonal image of $\mathbb{U}$.) If $V$ comes with a polarization form $Q$, then $Q(v \otimes v, v \otimes v) = Q(v,v)^2$ for the induced polarization on $V \otimes V$. Likewise, $\overline{Q}(v \otimes v) = \overline{Q}(v)^2$. Thus to bound $\NL_{q}$, we may reduce to bounding the Noether-Lefschetz locus $\NL^{\otimes 2}_{1,q}$ associated to $\mathbb{V} \otimes \mathbb{V}$ (since $\overline{Q}(v \otimes v) = \nu(v \otimes v)^2 = q^2$, as there is no square-free part). It then follows from \autoref{NLlocusbound} that
\[ \deg \NL_{q} \lesssim \deg \NL^{\otimes 2}_{1,q} \lesssim c \cdot q^{m^2 + \textrm{min}\{ r^2 + s^2, 2rs \} - 2} \]
where we use that the signature of the lattices appearing in the fibres of $\mathbb{V} \otimes \mathbb{V}$ is $(2rs, r^2 + s^2)$.

\section{Applications}

\subsection{Specialization of the Lower Bounds}
\label{specializationsec}

We would like to deduce \autoref{NLlocidegthmlower} and \autoref{abidemthm} from \autoref{mainthmlowbound} by taking the $T$ in the former two statements to be equal to the $S$ in the latter. However the statements of \autoref{NLlocidegthmlower} and \autoref{abidemthm} apply to \emph{any} sufficiently general hyperplane section in the sense of \autoref{suffgenhypplanedef}, whereas \autoref{mainthmlowbound} just works for some choice of such a section. Let us explain why in the case of our examples, this is not an issue, by following the proof of \autoref{mainthmlowbound} with $T = S$. We start by observing that, in the proof of \autoref{mainthmlowbound}, $L$ was required to be sufficiently general intersecting some $B' \subset B$ constructed in the proof. This $B'$ could be chosen to be any open subset of $B$ such that $g D_{M} \cap \psi(B')$ always has the ``expected'' dimension $d$ if non-empty. 

We observe that $B'$ has this property outside of a finite union of tensorial Hodge loci of $\mathbb{V}$. This follows from the ``Ax-Schanuel in families'' statement in \cite[\S3.3]{baldi2024effectiveatypicalintersectionsapplications}, where the bundle $P$ is constructed from $\mathbb{V}$ as in \cite[\S4.6]{baldi2024effectiveatypicalintersectionsapplications}. Indeed, the varieties $g \ch{D}_{M}$ belong to finitely many families of algebraic subvarieties of $\ch{D}$, hence induce finitely many algebraic families $f : \mathcal{Z} \to \mathcal{Y}$ of subvarieties of $P$ using the procedure outlined in \cite[\S4.6]{baldi2024effectiveatypicalintersectionsapplications}. Then the images in $T$ of intersections of fibres of $f$ with leaves of $P$, which in particular include all the germs of the form $\psi^{-1}(g D_{M})$, are contained in the fibres of finitely many families $\{ h_{i} : C_{i} \to A_{i} \}_{i=1}^{m}$ of weakly special subvarieties of $T$ as a consequence of \cite[Thm 3.12]{baldi2024effectiveatypicalintersectionsapplications}. We may assume that the fibres of the $h_{i}$ all have the same algebraic monodromy groups using \cite[Lem. 7.6]{baldi2024effectiveatypicalintersectionsapplications} (c.f. \cite[Lem. 7.5]{baldi2024effectiveatypicalintersectionsapplications}). The fibres of the $h_{i}$ which contain the $\psi^{-1}(g D_{M})$ do not have zero period dimension (since $\dim \varphi(\psi^{-1}(g D_{M})) > 0$ as $\dim D - \dim D_{M} < \dim \varphi(T)$ in both cases). It then follows from the $\mathbb{Q}$-simplicity of the algebraic monodromy group $\mathbf{H}_{T}$ in our two examples that each fibre of each $h_{i}$ lies in a weakly non-factor weakly special subvariety of $T$ in the sense of \cite[Def. 1.13]{fieldsofdef}. Then, as a consequence of \cite[Lem. 2.5]{fieldsofdef}, each such factor lies in a strict special subvariety of $T$. Thus the images of the $C_{i}$ in $T$ are contained in a union of strict special subvarieties (i.e., tensorial Hodge locus components) of $T$. There are countably many such varieties and the images of the $C_{i}$ are algebraic, so we can take this union to be finite.

\vspace{0.5em}

Now to prove \autoref{NLlocidegthmlower} and \autoref{abidemthm}, it suffices to replace $T$ with the complement in $T$ of this collection of tensorial Hodge loci. Note that by construction, any sufficiently general hyperplane section of the original $T$ also intersects this complement. After doing this, the locus $\mathcal{C}$ appearing in the proof of \autoref{mainthmlowbound} is automatically empty, and one can take $U$ to be the image of $q^{-1}(\psi(B))$ in $O(\mathbb{R})^{\circ}$. The statements of \autoref{NLlocidegthmlower} and \autoref{abidemthm} are thus reduced to an explicit estimate of the left-hand side of (\ref{ptcountasymp}), which we turn to in the remaining sections.

\vspace{0.5em}

Finally, let us note that the ``strongly $\mathbb{V}$-likely'' hypothesis in \autoref{mainthmlowbound} is checked in \cite[\S1.2.2]{arXiv:2303.16179} for the Noether-Lefschetz locus case, and is immediate from a simple dimension count in the split Jacobian case. The hypotheses on the normalizers is an easy consequence of \autoref{centrem} and the fact that the centralizers of the subgroups in question are finite.

\subsection{The Noether-Lefschetz Locus}

In this section we set $(a,b) = (r,s)$. We may identify $\mathbf{G}_{S} = \SO(a,b)$ and the natural inclusion $\mathbf{M} \hookrightarrow \mathbf{G}_{S}$ with $\mathbf{M} = \SO(a,b-1) \hookrightarrow \SO(a,b)$. 

We let $D_{v} \subset D_{S}$ be the natural Mumford-Tate subdomain consisting of Hodge structures with Mumford-Tate group $\SO(a,b-1)$ (the subscript $v$ denotes a primitive $\mathbb{Q}$-vector which $\SO(a,b-1)$ stabilizes). The upper bound appearing in \autoref{NLlocidegthmupper} is then a formal application of \autoref{NLlocusbound} together with the following lemma.

\begin{lem}
A Hodge subdatum $(\mathbf{M}, D_{M}) \subset (\SO(a,b), D_{S})$ which is conjugation-equivalent to the Hodge subdatum $(\SO(a,b-1), D_{v})$ is $\mathbb{Q}$-coset equivalent to the same subdatum if and only if $u(v) = u(v')$, where $v'$ is a non-zero element in the $\mathbb{Q}$-vector space stabilized by $\mathbf{M}$.
\end{lem}

\begin{proof}
By scaling $v$ we can assume that $Q(v,v) = u(v)$. We also let $g \in \SO(a,b)(\mathbb{R})$ be a point such that $(g^{-1} \mathbf{M} g, g^{-1} D_{M}) = (\SO(a,b-1), D_{v})$. 

Suppose first that $(\mathbf{M}, D_{M})$ is $\mathbb{Q}$-coset equivalent to $(\SO(a,b-1), D_{v})$. Then we may choose $g$ such that $v' = g \cdot v$ is defined over $\mathbb{Q}$. Choose $\lambda \in \mathbb{Q}^{\times}$ such that $\lambda^2 Q(v',v') = Q(\lambda v', \lambda v') = u(v')$. Then since $Q(v,v) = Q(g v, g v) = Q(v', v')$, necessarily $\lambda^2 Q(v,v) = u(v')$, and thus $\lambda^2 u(v) = u(v')$. Since $u(v)$ and $u(v')$ are both square free, it follows that $u(v) = u(v')$. 

Conversely, suppose that $u(v) = u(v')$. Because the subspace stabilized by $\SO(a,b-1)$ is $1$-dimensional, we have $g \cdot v = \lambda \cdot v'$ for some $\lambda \in \mathbb{R}^{\times}$, and we may assume $v'$ is chosen such that $Q(v',v') = u(v')$. We then have $Q(v,v) = Q(v',v')$. Then 
\[ Q(v,v) = Q(g v, g v) = \lambda^2 Q(v', v') = \lambda^2 Q(v,v) \]
so it follows that $\lambda = \pm 1$ and $g \cdot v$ is a $\mathbb{Q}$-vector. It follows that the coset $g \SO(a,b-1)$ is defined over $\mathbb{Q}$, or that the two Hodge data are $\mathbb{Q}$-coset equivalent. 
\end{proof}

We are therefore reduced to estimating the asymptotic number of rational points in a compact subset of $(\SO(a,b)/\SO(a,b-1))(\mathbb{R})$. We set $G = \SO(a,b)$ and $H = \SO(a,b-1)$ to match the notation of \S\ref{paramstripesec}. We start by giving an alternative description of the quotient $\SO(a,b)/\SO(a,b-1)$. First, consider the universal cover $\pi : \widetilde{G} \to \SO(a,b)$, which is an isogeny of $\mathbb{Q}$-algebraic groups. If we consider the scalar extension $\pi_{\mathbb{C}}$ we have that $\SO(a,b)_{\mathbb{C}} \cong \SO(a+b)_{\mathbb{C}}$. The fundamental group of $\SO(a+b)_{\mathbb{C}}$ is $\mathbb{Z}/2\mathbb{Z}$ and its universal cover is the Spin group $\textrm{Spin}(a+b)$; it follows that $\widetilde{G}$ is a $\mathbb{Q}$-form of $\textrm{Spin}(a+b)$.

Now consider $\widetilde{H} = \pi^{-1}(\SO(a,b-1))$. We claim that $\widetilde{H}$ is (geometrically) connected. Once again it suffices to scalar extend to $\mathbb{C}$, after which we can reduce to the same problem for $\SO(a+b-1) \subset \SO(a + b)$ and the universal covering by $\textrm{Spin}(a+b)$. Then the result follows from the fact that the construction of the Spin group is functorial in maps of quadratic spaces, which shows that the inverse image of $\SO(a+b-1)$ in $\textrm{Spin}(a+b)$ is identified with a copy of $\textrm{Spin}(a+b-1)$. 

We thus have an identification $\SO(a,b)/\SO(a,b-1) \cong \widetilde{G}/\widetilde{H}$. The groups $\widetilde{G}$ and $\widetilde{H}$ are both semisimple and simply-connected, so we may estimate the number of rational points a definable (relatively) compact subset of $(\SO(a,b)/\SO(a,b-1))(\mathbb{R})$ using the results of \S\ref{pointcountingcosetsec}. In particular, letting $\Omega$ be such a subset and applying our discussion in \S\ref{orbitasympsec} to the pair $(\widetilde{G}, \widetilde{H})$ and the natural representation $\widetilde{G} \to \GL_{m}$, we obtain from \autoref{pointcountingininterior} that 
\begin{equation}
\label{rationalpointcount}
(G \cdot v)(\mathbb{Q})^{\circ} \cap \Omega \cap \frac{1}{\nu} \mathbb{Z}^{m} \simeq \mu(B_{\nu}) = \mu(\Omega) \prod_{p} \mu_{p}(B_{\nu,p})
\end{equation}
where $B_{\nu,p} = \{ x_{p} \in G(\mathbb{Q}_{p}) \cdot v : \| x_{p} \|_{p} \leq p^{e_{p}} \}$. It therefore suffices to compute the volumes $\mu_{p}(B_{\nu,p})$. 

\subsubsection{Volume at $p$}
\label{volcompsubsubsec}

We recall that the variety $O = G \cdot v$ is a closed subvariety of $\mathbb{A}^{m}$ defined by the relation $Q(w,w) = Q(v,v) =: u$ for $w = (w_{1}, \hdots, w_{m})$ coordinates of $\mathbb{A}^{m}$. To compute using the measure $\mu_{p}$ on $O(\mathbb{Q}_{p})$ we may use Gauge forms. The Gauge forms $\Omega$ (resp $\omega)$ for the measures on $\SO(a,b)$ (resp. $O$) are given on \cite[pg. 149]{zbMATH00473017}. The expression for $\omega$ in particular, which is what interests us, is
\begin{align*}
\omega = \frac{dx_1 \wedge \cdots \wedge dx_{m}}{d(\overline{x}^{t} Q \overline{x})} .
\end{align*}
(Here we interpret the division of differential forms $\alpha/\beta$ as representing a differential form $\gamma$, if it exists, for which $\alpha = \gamma \wedge \beta$. The wedge product is computed in the algebra of differentials on $\mathbb{A}^{m}$.)

We write $L = \mathbb{Z}^{m}$ and let $B : L \to \mathbb{Z}$ be the quadratic form induced by $Q$. The lattice $(L, B)$ is unimodular and indefinite. Unimodular and indefinite lattices are fully classified by the signature of the form $Q$, and whether the lattice is ``even'' or ``odd'' (see \cite[Ch. II, Thm. 5.3]{zbMATH03456979}). (The lattice $(L, B)$ is even if $B(v,v) \in 2\mathbb{Z}$ for all $v, v' \in L$, and odd otherwise.) If the lattice is odd, then $Q$ can be diagonalized \cite[Ch. II, Thm. 4.3]{zbMATH03456979}. If we tensor with $\mathbb{Z}[2^{-1}]$ the distinction disappears, and $(L_{\mathbb{Z}[2^{-1}]}, B_{\mathbb{Z}[2^{-1}]})$ becomes isomorphic over $\mathbb{Z}[2^{-1}]$ to an odd unimodular indefinite lattice. In particular $(L, Q)_{\mathbb{Z}[2^{-1}]}$ is diagonalizable. As we are only interested in computing the factors in (\ref{rationalpointcount}) away from the prime $2$, we may therefore work over $\mathbb{Z}[2^{-1}]$ and assume that $Q$ is diagonal.

After diagonalizing $Q$ we may write $B(\overline{x}) = \overline{x}^t Q \overline{x} = \sum_{i} a_{i} x_{i}^2$. Then
\begin{align*}
d\left(\sum_{i} a_{i} x_{i}^2 \right) \wedge \frac{dx_2 \wedge \cdots \wedge dx_{m}}{x_1} &= 2 a_{1} \, x_{1} dx_{1} \wedge \frac{dx_2 \wedge \cdots \wedge dx_{m}}{x_1} \\
&= 2 a_{1} \, dx_{1} \wedge \cdots \wedge dx_{m} ,
\end{align*}
so we in fact have
\begin{align*}
\omega = \frac{1}{2 a_{1}} \frac{dx_2 \wedge \cdots \wedge dx_{m}}{x_1} .
\end{align*}

To integrate this form over the $p$-adics with $p \neq 2$ we apply the main result of \cite{zbMATH03853213}, which says that 

\begin{thm}
Let $f_{i} \in \mathbb{Z}_{p}[z_{1}, \hdots, z_{n}]$ be polynomials for $1 \leq i \leq N$ and write $\overline{f} : \mathbb{A}^{n} \to \mathbb{A}^{N}$ for the associated map. Let $t = (t_{1}, \hdots, t_{N}) \in \mathbb{Z}^{N}_{p}$ be a point such that the fibre $\overline{f}^{-1}(t)$ has at least one $\mathbb{Q}_{p}$-point at which the map $\overline{f}^{-1}(t) \to \mathbb{A}^{N}$ is submersive. Then one has
\begin{equation}
\int_{\overline{f}^{-1}(t) \cap \mathbb{Z}^{n}_{p}} | \Theta_{t} | = \lim_{e \to \infty} \frac{\# \{ a \in \mathbb{Z}^{n}_{p}/p^{e} : f_{i}(a) = t_{i} \textrm{ mod }p^{e} \textrm{ for }1 \leq i \leq N \} }{p^{(n-N)e}} =: \alpha_{p}(t, \overline{f}) ,
\end{equation}
where $\Theta_{t}$ is defined by 
\begin{equation}
\Theta_{t} = \textrm{sgn}(\lambda_{1}, \hdots, \lambda_{n})  \left| \det \left[ \frac{\partial f_{i}}{\partial z_{\Lambda}} \right] \right|^{-1} dz_{\Lambda^{c}}  .
\end{equation}
Here $\Lambda = (\lambda_{1}, \hdots, \Lambda_{N})$ is an ordered $N$-tuple of distinct increasing integers, and $\Lambda^{c} = (\lambda_{N+1}, \hdots, \lambda_{n})$ is the complementary tuple (c.f. \cite[\S2]{zbMATH03853213}). 
\end{thm}

We will apply this theorem with $N = 1$, $\Lambda = (1)$, $t = p^{2k} u$, $n = m$, and $f_{1} = B$ for some fixed integer $k \geq 0$. We consider the sets $A = B^{-1}(u) \cap \left(\frac{1}{p^{k}} \mathbb{Z}^{m}_{p}\right)$ and $A_{k} = B^{-1}(p^{2k} u) \cap \mathbb{Z}^{m}_{p}$. Let $z_{i} = p^{k} x_{i}$ for each $i$. We have that 
\begin{align*}
\int_{A} |\omega| &= \int_{A} \left| \frac{dx_{2} \wedge \cdots \wedge dx_{m}}{2 a_{1} x_1} \right| \\
&= \int_{A_{k}} \left| \frac{d(p^{-k} z_{2}) \wedge \cdots \wedge d(p^{-k} z_{m})}{p^{-k} 2 a_{1} z_{1}} \right| \\
&= p^{(m-2)k} \int_{A_{k}} \left| \frac{d z_{2} \wedge \cdots \wedge dz_{m}}{2 a_{1} z_{1}} \right| \\
&= p^{(m-2)k} \int_{A_{k}} \left| \Theta_{t} \right| \\
&= p^{(m-2)k} \, \alpha_{p}(p^{2k} u, B) 
\end{align*}

We now use the main theorem of \cite{zbMATH01486248} to compute $\alpha_{p}(p^{2k} u, B)$.

\begin{prop}
\label{alphaestprop}
We have
\begin{align*}
\alpha_{p}(p^{2k} u, B) &= 1 + \sum_{\substack{j=1 \\ j\textrm{ even}}}^{2k+1} p^{w(j)} + \sum_{\substack{j=1 \\ j\textrm{ odd}}}^{2k+1} p^{w(j)} \left(\frac{(- 1)^{m} \det Q}{p} \right) (\star) \\ (\star) &= \prod_{k = 1}^{m} \begin{cases} 
0 & 2k - j \geq 0, m \textrm{ odd} \\
(1-p^{-1}) \cdot \left(\frac{-1}{p}\right)^{m/2} & 2k - j \geq 0, m \textrm{ even} \\
\left(\frac{-1}{p}\right)^{(m+1)/2} \left(\frac{u}{p}\right) & 2k - j = -1, m \textrm{ odd} \\
- p^{-1/2} \cdot \left( \frac{-1}{p} \right)^{m/2} & 2k - j = -1, m \textrm{ even} ,
\end{cases} \\
w(j) &:= -mj/2 + j + \textrm{min}\{ 2k-j,0 \}/2 .
\end{align*}
\end{prop}

\begin{proof}
This is just a specialization of the main theorem of \cite{zbMATH01486248}; since the expression appearing there is complicated, we give some details, following the notation appearing there. We thus set $T = p^{2k} u$ and $S = Q$. We write 
\[ \alpha_{p}(T, S) = \sum_{?} F_{1}(?) \sum_{k=0}^{r+1} F_{2}(?, k) \sum'_{\{ \nu \}_{k}} F_{3}(?, k, \{ \nu \}_{k}) \]
in accordance with the expression for $\alpha_{p}(T,S)$ appearing in \cite{zbMATH01486248}. We have $n = 1$, so $\mathfrak{S}_{n} = 1$ is the trivial group; $c_{1}(\sigma) = 1$ always and $c_{2}(\sigma) = 0$ always; the set $I = \{ 1 \}$ and admits a unique $\mathfrak{S}_{n}$ stable partition of length $r+1 = 1$. This means that the functions $t$ and $\tau$ are always zero. From all this it follows that $F_{1} = 1/2$, the first summation sign can be removed, and we are left to estimate an expression of the form
\[ \alpha_{p}(T,S) = \frac{1}{2} \sum_{k=0}^{1} F_{2}(k) \sum'_{\{ \nu \}_{k}} F_{3}(k, \{ \nu \}_{k}) . \]

Now $n_{\ell}, n^{(\ell)}$ and $n(\ell)$ are equal $1$ when $\ell \leq 0$, and when $\ell \geq 1$ we have $n_{\ell} = n^{(\ell)} = n(\ell) = 0$. Moreover $c^{(0)}_{1}$ is always just $c_{1}$, so equal to $1$, and $c^{(1)}_{1} = 0$. It follows that $F_{2}(0) = 2$ and $F_{2}(1) = 1$. When $k = 0$ the third summation is equal to $1$ by convention, so we obtain
\[ \alpha_{p}(T,S) = 1 + \frac{1}{2} \sum'_{\{ \nu \}_{1}} F_{3}(1, \{ \nu \}_{1}) . \]
Now $\{ \nu \}_{1}$ is just a single integer $\nu$ satisfying $- b_{0}(\textrm{id}, T) \leq \nu \leq -1$, where $b_{0}(\textrm{id}, T) = \beta_{1} + 1$. We thus get
\begin{align*}
\sum'_{\{ \nu \}_{1}} F_{3}(1, \{ \nu \}_{1}) &= \sum_{-\nu = 1}^{\beta_{1} + 1} p^{-\nu} \Xi_{0,\nu}(\textrm{id}; T, S) .
\end{align*}
Now with $\nu = \lambda < 0$ one has $\rho_{0,\lambda} = \frac{m \lambda}{2} + \frac{1}{2} \textrm{min} \{ \beta_{1} + \lambda, 0 \}$, and $\Xi_{0, \lambda}(\textrm{id}, T, S) = p^{\rho_{0,\lambda}} \xi_{1,\lambda}(T,S)$. In our setting $B_{i}(\lambda) = \varnothing$, and $A(\lambda) = \{ k : 1 \leq k \leq m, \lambda \neq 0 \textrm{ mod }2 \}$, i.e., either $A(\lambda) = \{ 1, \hdots, m \}$ or $A(\lambda) = \varnothing$, depending on the parity of $\lambda$. 

Thus $\xi_{1, \lambda}(T,S)$ is equal to $2$ when $\lambda$ is even, and for $\lambda$ odd is given by the expression
\begin{align*}
\xi_{1,\lambda} = 2 \left(\frac{(-1)^{m} \det Q}{p} \right) \prod_{k = 1}^{m} \begin{cases} 
0 & \beta_{1} + \lambda \geq 0, m \textrm{ odd} \\
(1-p^{-1}) \cdot \left(\frac{-1}{p}\right)^{m/2} & \beta_{1} + \lambda \geq 0, m \textrm{ even} \\
\left(\frac{-1}{p}\right)^{(m+1)/2} \left(\frac{u}{p}\right) & \beta_{1} + \lambda = -1, m \textrm{ odd} \\
- p^{-1/2} \cdot \left( \frac{-1}{p} \right)^{m/2} & \beta_{1} + \lambda = -1, m \textrm{ even} .
\end{cases}
\end{align*}
Putting everything together we obtain the result in the statement.
\end{proof}

\begin{cor}
We have $\alpha_{p}(p^{2k} u, B) \geq \frac{1 - 2 p^{(1-m/2)}}{1 - p^{1-m/2}}$.
\end{cor}

\begin{proof}
From \autoref{alphaestprop} we get that
\begin{align*}
\alpha_{p}(p^{2k} u, B) &\geq 1 - \sum_{j=1}^{2k+1} p^{(1-m/2) j} \\
&\geq 1 - \left(\frac{1}{1 - p^{(1-m/2)}} - 1 \right) \\
&= \frac{1 - 2 p^{(1-m/2)}}{1 - p^{1-m/2}} .
\end{align*}
\end{proof}

As a consequence of the corollary and the calculation preceding \autoref{alphaestprop}, one has $\mu_{p}(B_{\nu,p}) \geq \left( \frac{1 - 2 p^{1 - m/2}}{1 - p^{1-m/2}} \right) p^{(m-2) k}$. Taking the product over all $p$ the product $\prod_{p} \left( \frac{1 - 2 p^{1 - m/2}}{1 - p^{1-m/2}} \right)$ is just a positive constant, so one obtains $\prod_{p} \mu_{p}(B_{\nu,p}) \geq c \cdot \nu^{m-2}$ for some constant $c > 0$, hence the result.

\subsection{Splittings of Jacobians}

In this case we have $(\mathbf{G}_{S} = \Sp_{g}, \mathbb{H}_{g})$, with $\mathbb{H}_{g}$ the Siegel upper half space, and $(\mathbf{M} = \Sp_{k} \times \Sp_{g-k}, \mathbb{H}_{k} \times \mathbb{H}_{g-k})$. 

\begin{lem}
There are only finitely many $\mathbb{Q}$-coset equivalence classes of Hodge subdata of $(\Sp_{g}, \mathbb{H}_{g})$ which are isomorphic to $(\Sp_{k} \times \Sp_{g-k}, \mathbb{H}_{k} \times \mathbb{H}_{g-k})$. 
\end{lem}

\begin{proof}
Because there are finitely many equivalence classes of Hodge subdata under real-conjugation equivalence, it suffices to fix such an equivalence class and show that it contains finitely many sub-equivalence classes under $\mathbb{Q}$-coset equivalence. Given two such subdatum $(\mathbf{M}, D_{M})$ and $(\mathbf{M}', D_{M'})$ of the Hodge datum $(\Sp_{g}, \mathbb{H}_{g})$, one knows that the subset of $\textrm{End}(\mathbb{Q}^{2g})$ stabilized by the conjugation action of $\mathbf{M}$ (resp. $\mathbf{M}'$) is the $2$-dimensional subspace spanned by $\{ e, 1 - e \}$ (resp. $\{ e', 1-e' \}$) where $e : \mathbb{Q}^{m} \to \mathbb{Q}^{m}$ (resp. $e' : \mathbb{Q}^{m} \to \mathbb{Q}^{m}$) is a non-trivial idempotent. A real element $g \in \Sp_{g}(\mathbb{R})$ that conjugates the two Hodge subdata then induces an equality of algebras 
\[ g\, \textrm{span}_{\mathbb{R}}\{ e, 1-e \} g^{-1} = \textrm{span}_{\mathbb{R}}\{ e', 1-e' \} . \]
In particular $g e g^{-1}$ is a non-trivial idempotent in the algebra $\textrm{span}_{\mathbb{R}}\{ e', 1-e' \}$, necessarily equal to either $e'$ or $1-e'$. Since $\mathbf{M}$ (resp. $\mathbf{M}'$) is exactly the stabilizer in $\Sp_{g}$ of $e$ (resp. $e'$), it follows after applying the isomorphism $\Sp_{g}/\mathbf{M} \simeq \Sp_{g} \cdot e$ that the coset $g \mathbf{M}$ is defined over $\mathbb{Q}$. 
\end{proof}

To deduce our results, we compare the degrees of isogenies associated to $k$-factors of an abelian variety $A$ with denominators of Hodge-theoretic idempotents. 

\begin{defn}
For a weight one integral Hodge structure $V$, we say that $V$ has a $k$-factor of denominator $\nu$ if there exists a rank $2k$-Hodge-theoretic idempotent $e : V_{\mathbb{Q}} \to V_{\mathbb{Q}}$ such that $\nu e \in \textrm{End}(V)$ is primitive. 
\end{defn}

\begin{lem}
\label{isoheightisdenom}
Suppose that $A$ is a complex abelian variety, and that $V$ is the associated weight one integral Hodge structure. Then 
\begin{itemize}
\item[-] if $A$ has a $k$-factor of degree dividing $\nu$, then $V$ has a $k$-factor of denominator dividing $\nu$; and
\item[-] if $V$ has a $k$-factor of denominator dividing $\nu$, then $A$ has a $k$-factor of degree dividing $\nu^{2g}$.
\end{itemize}
\end{lem}

\begin{proof}
Suppose first that $A$ has a $k$-factor of degree $d$ dividing $\nu$. Let $V'$ be an integral weight one Hodge structure with $\alpha : V \to V'$ an isogeny of degree $d$ such that $V' = V_{1} \oplus V_{2}$ is a decomposition of weight $1$ integral Hodge structures, and $V_{1}$ has rank $2k$. Let $e' : V' \to V_{1}$ be the associated Hodge-theoretic idempotent. Then $e = \alpha^{-1} e' \alpha$ is a Hodge-theoretic idempotent $V_{\mathbb{Q}} \to V_{\mathbb{Q}}$ of rank $2k$ and $(\det \alpha) e$ is integral, so the denominator of $e$ divides $(\det \alpha) = d$, and $d | \nu$ by assumption.

For the second statement we start with a rank $2k$ Hodge-theoretic idempotent $e : V_{\mathbb{Q}} \to V_{\mathbb{Q}}$ of $V$ such that $\nu e \in \textrm{End}(V)$ and construct an isogenous Hodge structure $V'$. As a lattice, we define $V' = \textrm{span}_{\mathbb{Z}} \{ e(V), (1-e)(V) \}$. It is clear that $V'$ is preserved by $e$, and that the Hodge decomposition on $V_{\mathbb{C}}$ induces a Hodge-decomposition on $V'_{\mathbb{C}}$. Moreover one has $V \subset V'$, since $v = e v + (1-e) v$ for any $v \in V$. Then the inclusion $V \subset V'$ is an isogeny $\alpha$. Its degree is the index $[V' : V]$. We then have
\[ \nu^{2g} = \left[ \frac{1}{\nu} V : V \right] = \left[ \frac{1}{\nu} V : V' \right] [V' : V] , \]
hence $(\deg \alpha) | \nu^{2g}$. 
\end{proof}

We define 
\begin{align*}
\textrm{SP}'_{k,\nu} := \{ s \in S(\mathbb{C}) : \mathbb{V}_{s} \textrm{ has a }k\textrm{-factor of denominator dividing }\nu \} . 
\end{align*}
Note that our lemma shows that $\textrm{SP}_{k,\nu} \subset \textrm{SP}'_{k,\nu}$ and $\textrm{SP}'_{k,\nu} \subset \textrm{SP}_{k,\nu^{2g}}$. 

Our result \autoref{abidemthmupper} now follows immediately from \autoref{NLlocusbound} and \autoref{isoheightisdenom}, where we note that components of the split Jacobian locus become components of the Noether-Lefschetz locus upon replacing $\mathbb{V}$ with $\textrm{End}(\mathbb{V})$. In particular, one has 
\begin{align*}
2 \dim O_{k} &= 2 ( \dim \Sp_{g} - \dim (\Sp_{k} \times \Sp_{g-k}) )   \\
&= 2 ( g(2g+1) - (k(2k+1) + (g-k)(2(g-k)+1)) ) \\
&= 8k(g-k) .
\end{align*}

For the lower bound, we may first prove it for $\textrm{SP}'_{k,\nu}$, and then use that $\textrm{SP}_{k,\nu^{2g}}$. For this we use that both $\mathbf{G}_{S}$ and $\mathbf{M}$ are semisimple and simply-connected, so we can reason exactly as before to deduce the result from \autoref{pointcountingininterior} and \autoref{mainthmlowbound}, our discussion in \S\ref{specializationsec}.

\bibliography{hodge_theory}
\bibliographystyle{alpha}

\end{document}